\documentclass[12pt,a4paper]{amsart}
\usepackage[latin1]{inputenc}
\usepackage[english]{babel}
\usepackage{amsfonts,amsmath}
\usepackage{mathabx}
\usepackage{fancyhdr,graphicx, xcolor}
\usepackage{hyphenat}
\usepackage{amsthm} 

\usepackage[scanall,209mode]{psfrag}

\def\R{\mathbb R}
\def\C{\mathbb C}
\def\N{\mathbb{N}}

\def\D{\text{Diag}}
\def\col{\text{col}}

\DeclareMathOperator{\Diag}{Diag}

\def\tr{\text{tr}}

\def\restrict#1{\raise-.5ex\hbox{\ensuremath{\big|}}_{#1}} 

\usepackage[hidelinks]{hyperref} 
\newtheorem{theorem}{Theorem}

\newtheorem{corollary}{Corollary}

\newtheorem{definition}{Definition}
\newtheorem{example}{Example}

\newtheorem{lemma}{Lemma}

\newtheorem{proposition}{Proposition}
\newtheorem{remark}{Remark}

\numberwithin{equation}{section}

\begin{document} 
	\title[Minimal compact operators, subdifferential and SDP ]{Minimal compact operators, subdifferential of the maximum eigenvalue and semi-definite programming}
	
	\author{Tamara Bottazzi $^{1,2}$ and Alejandro Varela$^{3,4}$}

	\address{$^1$ Universidad Nacional de R\'io Negro. Centro Interdisciplinario de Telecomunicaciones, Electrónica, Computación y Ciencia Aplicada, Sede Andina (8400) S.C. de Bariloche, Argentina.}
	\address{$^2$ Consejo Nacional de Investigaciones Cient\'ificas y T\'ecnicas, (1425) Buenos Aires,
		Argentina.}
	
	\address{$^3$Instituto Argentino de Matem\'atica ``Alberto P. Calder\'on", Saavedra 15 3er. piso, (C1083ACA) Buenos Aires, Argentina}
	
	\address{$^4$Instituto de Ciencias, Universidad Nacional de Gral. Sarmiento, J.	
		M. Gutierrez 1150, (B1613GSX) Los Polvorines, Argentina}
	
	\email{tbottazzi@unrn.edu.ar, avarela@campus.ungs.edu.ar }
	
	\thanks{Partially supported by Grants CONICET (PIP 0525), ANPCyT (PICT 2015-1505 and 2017-0019)}
	
	\maketitle

	\begin{abstract}
		We formulate the issue of minimality of self-adjoint operators on a Hilbert space as a semi-definite problem, linking the work by Overton in \cite{overton} to the characterization of minimal hermitian matrices. This motivates us to investigate the relationship between minimal self-adjoint operators and the subdifferential of the maximum eigenvalue, initially for matrices and subsequently for compact operators. In order to do it we obtain new formulas of subdifferentials of maximum eigenvalues of compact operators that become useful in these optimization problems.
		
		Additionally, we provide formulas for the minimizing diagonals of rank one self-adjoint operators, a result that might be applied for numerical large-scale eigenvalue optimization.
%
	\end{abstract}

\section{Introduction and preliminaries}
Let $B(H)$ and $K(H)$ be the spaces of linear bounded and compact operators defined on a Hilbert space $H$, respectively. We call $A\in B(H)$ a minimal operator if $\|A\|\leq\|A+D\|$, for all $D$ diagonal in a fixed orthonormal basis $E=\{e_i\}_{i\in I}$ of $H$ and $\|\cdot\|$ the operator norm. Note that when $A\in K(H)$, we can suppose that $H$ is separable since there is only a numerable set $\{e_{i_k}\}_{k\in \N}$ such that $A(e_{i_k})\neq 0$. 

In particular, minimal compact self-adjoint operators are related with the distance to the subspace of diagonal self-adjoint operators, denoted by $\D(K(H))^{sa}$, since for $A\in K(H)$ self-adjoint
$$
\text{dist }(A,\D(K(H))^{sa})=\inf_{D\in \Diag(K(H))}\|A+D\|.
$$
Minimal operators allow the concrete description of geodesics in homogeneous spaces obtained as orbits of unitaries under a natural Finsler metric (see  \cite{dmr1}).

In the case that $H=\C^n$, $B(H)$ is the space of complex square matrices of $n\times n$, that is $M_n(\C)$. The matricial case of minimal operators was extensively studied in \cite{amlmrv}, \cite{alrv} and \cite{KV 3x3}. 

In \cite{bot var best}, \cite{bot var min length curves in unitary orb} and  \cite{BV Dga} we studied minimal self-adjoint compact operators where it was stated that in general neither existence nor uniqueness of compact minimizing diagonals was granted. Some of these results were recently generalized to more general subalgebras of $K(H)$ and to C$^*$-algebras in \cite{zhang-jiang, Zhang-Jiang-2023}.

The characterization of minimal self-adjoint matrices can be stated as a semi-definite programming problem \cite{overton}. Moreover, in \cite{overton2}  Overton develops several algorithms using the subdifferential of the maximum eigenvalue of a matrix, which is the set
$$\partial \lambda_{max}(A)=\left\lbrace V\in M_n(\C): V=V^*\ \text{and } \lambda_1(Y)-\lambda_1(A)\geq \text{Re } \tr(V(Y-A)), \ \forall\ Y\in M_n^h(\C)\right\rbrace .
$$
This subdifferential was also studied in \cite{watson} and is a powerful tool in cases of non-differentiable functions \cite{clarke1,clarke2}.

The work of Overton in \cite{overton} and \cite{overton2} motivated us to study the relation between minimal operators and subdifferentials, first for matrices, and then for compact operators. 

In \cite{bata-grover, Grover, singla}, the authors give useful expressions for the subdifferential of the norm operator and they relate this concept with the distance to some closed subsets in $B(H)$. 

In the present work, we relate minimal operators with subdifferentials of the maximum eigenvalue and of the norm. We vinculate these concepts with the moment of the eigenspace of the maximum eigenvalue and the joint numerical range, which was developed in \cite{KV momento} for matrices, and \cite{BV-Moment} for compact operators. 

Indeed, one of our main results is an explicit formulation for the largest eigenvalue subdifferential of a compact self-adjoint operator $A(x)$ with variable real diagonal $x$,  
\begin{equation*}
		\partial\big(\lambda_{max}(A(x))\big)=\Diag(\partial\lambda_{max}(A(x)))= m_{S_{max}},
\end{equation*}
in terms of $m_{S_{max}}=\text{co}\{|v|^2: v\in S_{max}, \|v\|=1\}$, the moment of the eigenspace related to the largest eigenvalue $\lambda_{max}(A(x))$ (see \eqref{defi moment Y matrix}, \eqref{defi moment Y compact} and  Theorem \ref{teo subdif=moment}). Additionally, when the smallest eigenvalue $\lambda_{min}(A(x))$ is negative, we give explicit formulas for $\partial\big(\lambda_{min}(A(x))\big)$ and vinculate $\partial\lambda_{max}$ and $\partial\lambda_{min}$ with the subdifferential of the spectral norm of $A(x)$.
The above leads us to a new characterization of minimal self-adjoint compact operators that involves $\partial\big(\lambda_1(A(x))\big)$, $\partial\|A(x)\|$, the intersection of moments of the maximum and minimum eigenvalues of $A(x)$ and the joint numerical range of a certain family of operators that has been studied in \cite{BV-Moment} and \cite{KV momento}.

We first obtain the subdifferential formulas and the characterization of minimal operators when $A(x)\in M^{h}_n(\C)$. In order to extend it to the compact operator case, we needed some additional tools from non-smooth analysis and optimization. We used \cite{clarke1} and \cite{clarke2} as our main references of the topic. 

The formula of the subdifferentials that we obtain can be applied to eigenvalues with multiplicity higher than one, but in case of a simple eigenvalue, our formula coincides with the definition of gradient and partial derivatives of the maximum eigenvalue of a matrix (see \cite{lancaster} and \cite{kangal}). This may be useful to develop or improve algorithms for large-scale eigenvalue optimization.

The results we present in this paper are divided in three parts. 
Section \ref{sec sdp} is devoted to state the minimality of self-adjoint matrices as a semi-definite problem, relating some of the results obtained by Overton in \cite{overton} with the main characterization theorems that appear in \cite{alrv}.
Inspired by \cite{overton2}, in Section \ref{sec subdif} we study the subdifferential of the maximum eigenvalue, first for matrices and then for compact operators, and we link it with the moment of its eigenspace. In order to obtain these results we calculate new formulas of subdifferentials of eigenvalues of compact self-adjoint operators and operator norms (see Theorems \ref{teo subdif y moment compact} and \ref{equiv subdif, momentos, etc para compactos}).
Finally, in section \ref{sec rank one} we show explicit formulas for the minimizing diagonals for a given rank-one self-adjoint compact operator. These results might be used to improve some algorithms recently obtained for large-scale eigenvalue optimization problems (see \cite{kangal}).

Next we introduce some additional definitions and notations.

We use the superscript $^{sa}$ to note the subset of self-adjoint elements of a particular subset of $B(H)$. A self-adjoint element $A\in B(H)$ is called positive if  $\langle Ax, x\rangle \geq 0$ for all $x\in H$ and it is denoted by $A\geq 0$. For an  operator $A \in B(H)$ we use $\ker(A)$ to denote the kernel of $A$ and $|A|$ the modulus of $A$ given by $(A^*A)^{1/2}$.

For every compact operator $A\in K(H)$, let $s_1(A), s_2(A),\cdots $ be the singular values of $A$, i.e.
the eigenvalues of $|A|$  in decreasing order ($s_i(A)=\lambda_i(|A|)$, for each $i\in \N$) and repeated
according to multiplicity. Let
\begin{equation} \label{defipllel}
	{\|A\|}_1 = \sum_{i = 1}^\infty s_i(A) = {\rm tr}|A|,
\end{equation}
where $\rm tr(\cdot)$ is the trace functional, i.e.
\begin{equation}\label{traza}
	{\rm tr}(A)=\sum_{j=1}^{\infty} \langle Ae_j,e_j\rangle
\end{equation}
where $e_j$ are the elements of a fixed orthonormal basis $E$.
Observe that the series \eqref{traza} converges
absolutely and it is independent from the choice of basis and this coincides with the usual definition of the trace if $H$ is finite-dimensional. 

We define the usual ideal of trace class operators as
\begin{equation}
B_1(H)=\{A\in K(H):\ \|A\|_1<\infty\}.
\end{equation}

\section{The characterization of minimal matrices as a semidefinite-problem} \label{sec sdp}
Let $A_0\in M_n^h(\C)$ and $\varphi:\R^n\to \R$ be the function given by
$$\varphi(x)=\max_{1\leq i\leq n}|\lambda_i(A(x))|,$$
where $A(x)= A_0+\D(x)$, $x\in \R^n$ and $\{\lambda_i(A(x))\}_{i=1}^n$ are its eigenvalues in decreasing order counted with multiplicity. We are interested to study the following convex optimization problem:
\begin{equation}\label{P1}
\min_{x\in \R^n}\varphi(x).
\end{equation}
Any solution of \eqref{P1} gives us the spectral norm of a minimal matrix $A(x_0)$ and a best real diagonal approximation to the subspace of real diagonal matrices $\D(x_0)$ (may be not unique). In this case, we say that $A(x_0)$ is a minimal matrix, that is
$$\|A(x_0)\|\leq \|A(x)\|,\ \text{for every } x\in \R^n.$$

When $A_0$ is a real symmetric matrix, this problem is a particular case from \cite{overton} and it can be stated as
\begin{equation}\label{P2}
\min_{w\in \R,\ x\in \R^n}w\text{ such that } -w\leq \lambda_i(A(x))\leq w,\ 1\leq i\leq n,
\end{equation}
or equivalently
\begin{equation}\label{P3}
\min_{w\in \R,\ x\in \R^n}w\text{ such that }\left\lbrace \begin{array}{ll}
 wI-A(x)\geq 0\\
 wI+A(x)\geq 0.
\end{array}\right. 
\end{equation}
Problem \eqref{P3} can be viewed as a Semi-Definite Programming (SDP) issue with two semidefinite constraints. Fletcher in \cite{fletcher} deals with a similar problem with only one semidefinite constraint. 

If $A(x)$ is minimal, then there exist natural numbers $1\leq s,t<n$ such that 
\begin{equation}\label{st}
\left\lbrace \begin{array}{ll}
	\lambda_i(A(x))=\lambda_i&i=1,2,...,n\\
	\lambda_i=w&i=1,2,...,t\\
	\lambda_i=-w&i=n-s+1,...,n\\
	w=\lambda_1=...=\lambda_t&>\lambda_{t+1}\geq...\geq \lambda_{n-s}>\lambda_{n-s+1}=...=\lambda_n=-w.
\end{array}\right. 
\end{equation}
If $\{q_1,...,q_n\}$ is an orthonormal set of eigenvectors corresponding to $\{\lambda_1,...,\lambda_n\}$, we define 
\begin{equation}
\label{def Q1 y Q2}
Q_1=[q_1|...|q_t] \text{ and } Q_2=[q_{n-s+1}|...|q_n], 
\end{equation}
matrices of $n\times t$ and $n\times s$, respectively.

Let $E_{k}=e_k\otimes e_k=e_ke_k^*$, with $\{e_k\}_{k=1}^{m}$ the canonical basis of $\R^m$ for any $m\in \mathbb{N}$.

The next result is a particular case of Theorem 3.2 in \cite{overton} applied to our context. 
\begin{proposition}
Let $A_0\in \R^{n\times n}$, $A_0=A_0^t$ and $x\in \R^n$. The following conditions are equivalent:
\begin{enumerate}
	\item $x$ is a solution of \eqref{P2} (i.e: $A(x)$ is minimal).
	\item $A(x)$ fulfills \eqref{st} and there exist semidefinite positive symmetric matrices $U$ of $t\times t$ and $V$ of $s\times s$ such that
	\begin{itemize}
		\item $\tr(U)+\tr(V)=1$,
		\item $\tr(Q_1^tE_kQ_1U)-\tr(Q_2^tE_kQ_2V)=0$, $k=1,...,n$.
	\end{itemize}
\end{enumerate}
\end{proposition}

According to \cite{boyd-book} we can convert an SDP complex problem into a real SDP, using the following result.
\begin{lemma} \label{pos}
	For every $Y\in  M_n^h(\C)$, 
	\begin{equation}\label{eqpos}
	Y\geq 0 \text{ if and only if } \begin{bmatrix}
	\Re(Y)&-\Im(Y)\\
	\Im(Y)&\Re(Y)
	\end{bmatrix}\geq 0,
	\end{equation}
	(here $\Re(Y)=\frac{1}{2}(Y+\overline{Y})$ and $\Im(Y)=\frac{1}{2i}(Y-\overline{Y})$).
	
\end{lemma}

\begin{proof}
First define the block matrix $U=\frac{1}{\sqrt{2}}\begin{bmatrix}
I_n&iI_n\\
iI_n&I_n
\end{bmatrix}$,
and observe that 
$$UU^*=U^*U=\begin{bmatrix}
I_n&0\\
0&I_n
\end{bmatrix}.$$
Then, 
$$U^*\begin{bmatrix}
Y&0\\
0&\overline{Y}
\end{bmatrix}U=\frac12\begin{bmatrix}
Y+\overline{Y}&i(Y-\overline{Y})\\
-i(Y-\overline{Y})&Y+\overline{Y}
\end{bmatrix}=\begin{bmatrix}
\Re(Y)&-\Im(Y)\\
\Im(Y)&\Re(Y)
\end{bmatrix}.$$
Therefore, $\begin{bmatrix}
\Re(Y)&-\Im(Y)\\
\Im(Y)&\Re(Y)
\end{bmatrix}$ is unitary block equivalent to $\begin{bmatrix}
Y&0\\
0&\overline{Y}
\end{bmatrix}$ and 
$$\begin{bmatrix}
\Re(Y)&-\Im(Y)\\
\Im(Y)&\Re(Y)
\end{bmatrix}\geq 0\ \text{ if and only if }\begin{bmatrix}
Y&0\\
0&\overline{Y}
\end{bmatrix}\geq 0.$$
So, in order to prove \eqref{eqpos}, we only need to show that 
$$Y\geq 0\ \text{ if and only if }\begin{bmatrix}
Y&0\\
0&\overline{Y}
\end{bmatrix}\geq 0.$$
It is evident that the inequality 
$$\begin{bmatrix}
Y&0\\
0&\overline{Y}
\end{bmatrix}\geq 0$$
implies $Y\geq 0$. On the other hand, if $Y\geq 0$, there exist a unitary $V\in M_n(\C)$ such that $Y=V^*\D(\lambda(Y))V$, with $\lambda_i(Y)\geq 0$, for every $i=1,...,n.$ Then, 
$$\overline{Y}=\overline{V}^*\D(\overline{\lambda(Y)})\overline{V}=\overline{V}^*\D(\lambda(Y))\overline{V}\geq 0,$$
since $\overline{V}$ is also unitary, and
$$\begin{bmatrix}
Y&0\\
0&\overline{Y}
\end{bmatrix}\geq 0.$$
\end{proof}

\begin{proposition} \label{sdp real of 2n}
	Let $A_0\in  M_n^h(\C)$, $x\in \R^n$ and $\bar{A}[x]\in \mathbb{R}^{2n\times 2n} $ such that
	\begin{equation}\label{A 2n x 2n}
	\bar{A}[x]=\begin{bmatrix}
	\Re(A_0)&-\Im(A_0)\\
	\Im(A_0)&\Re(A_0)
	\end{bmatrix} +\begin{bmatrix}
	\D(x)&0\\
	0&\D(x)
	\end{bmatrix}.
	\end{equation}
	Then, we state problem \eqref{P1} as the following real SDP
	\begin{equation}\label{P4}
		\min_{w\in \R,\ x\in \R^n}w\text{ such that }
	\end{equation}

	\begin{equation}\label{cond}
		w \begin{bmatrix}
			I&0\\
			0&I
		\end{bmatrix} - \bar{A}[x]\geq 0\ \text{and }\ 	w \begin{bmatrix}
		I&0\\
		0&I
		\end{bmatrix} + \bar{A}[x]\geq 0.
	\end{equation}
	
\end{proposition}
	
\begin{proof}
By Lemma \ref{pos}, linear matrix restrictions $wI- A(x)\geq 0$ and $wI+ A(x)\geq 0$ are equivalent to conditions in \eqref{cond}. 

\end{proof}

If $\bar{A}[x]$ is a solution of \eqref{P4} and \eqref{cond}, then there exist natural numbers $1\leq \bar{s},\bar{t}<2n$ such that 
\begin{equation}\label{st2n}
\left\lbrace \begin{array}{ll}
\lambda_i(\bar{A}[x])=\lambda_i&i=1,2,...,2n\\
\lambda_i=w&i=1,2,...,\bar{t}\\
\lambda_i=-w&i=n-\bar{s}+1,...,2n\\
w=\lambda_1=...=\lambda_{\bar{t}}&>\lambda_{\bar{t}+1}\geq...\geq \lambda_{2n-\bar{s}}>\lambda_{2n-\bar{s}+1}=...=\lambda_{2n}=-w.
\end{array}\right. 
\end{equation}
If $\{q_1,...,q_{2n}\}$ is an orthonormal set of eigenvectors corresponding to $\{\lambda_1,...,\lambda_{2n}\}$, we define 
\begin{equation}
\label{def Q1 y Q2 2n}
\bar{Q_1}=[q_1|...|q_{\bar{t}}] \text{ and } \bar{Q_2}=[q_{n-\bar{s}+1}|...|q_{2n}], ,
\end{equation}
matrices of $2n\times \bar{t}$ and $2n\times \bar{s}$, respectively.

Therefore, we arrive to the next result in relation with the study of
\begin{equation}\label{PBlock}
	\min_{x\in \R^n}\|\bar{A}[x]\|.
\end{equation}

\begin{theorem}[SDP Complex into SDP real]\label{complex to real}
	Let $A_0\in  M_n^h(\C)$, $x\in \R^n$ and $\bar{A}[x]\in \mathbb{C}^{2n\times 2n}$ as in \eqref{A 2n x 2n}. The following conditions are equivalent:
	\begin{enumerate}
		\item $x$ is a solution of \eqref{P1} (i.e. $A(x)$ is a minimal matrix of $n\times n$). 
		\item $\|\bar{A}[x]\|\leq \|\bar{A}[y]\|$ for every $y\in \R^n$ (i.e.,  $(x,x)\in \R^{2n}$ is a solution of \eqref{PBlock}).
		\item Following the same notation as in Proposition \ref{sdp real of 2n}, \eqref{st2n} and \eqref{def Q1 y Q2 2n}, there exist semidefinite positive symmetric matrices $\bar{U}$ of $\bar{t}\times \bar{t}$ and $\bar{V}$ of $\bar{s}\times \bar{s}$ such that
		\begin{itemize}
			\item $\tr(\bar{U})+\tr(\bar{V})=1$,
			\item $\tr(\bar{Q_1}^t(E_k+E_{n-k})\bar{Q_1}\bar{U})-\tr(\bar{Q_2}^t(E_k+E_{n-k})\bar{Q_2}\bar{V})=0$, $k=1,...,n$.
		\end{itemize}
		
	\end{enumerate}
\end{theorem}
\begin{proof}
It follows directly from the conversion of Problem \eqref{P1} into a SDP real problem, as we did in Proposition \ref{complex to real}.
\end{proof}

Theorem \ref{complex to real} indicates that $A(x)$ is minimal if and only if the real $2n\times2n$ block matrix $\bar{A}[x]$
is a solution of the problem \eqref{PBlock}.

Observe that a solution $\bar{A}[x]$ of \eqref{PBlock} is not necessarily a minimal matrix of $2n\times2n$, since it can exist a $2n\times 2n$ best real diagonal approximant $D$ such that 
$$ \|\bar{A}[0]+D\|<\|\bar{A}[x]\|\text{ with }
D\neq \begin{bmatrix}
	\D(x)&0\\
	0&\D(x)\end{bmatrix}=\sum_{k=1}^{n}x_kA_k,$$
where $A_k=E_k+E_{n-k}$ and
\begin{eqnarray*}
	\bar{A}[x]&=&\begin{bmatrix}
		\Re(A_0)&-\Im(A_0)\\
		\Im(A_0)&\Re(A_0)
	\end{bmatrix} +x\sum_{k=1}^{n}e_ke_k^t+x\sum_{k=n+1}^{2n}e_ke_k^t\\
	&=&\begin{bmatrix}
		\Re(A_0)&-\Im(A_0)\\
		\Im(A_0)&\Re(A_0)
	\end{bmatrix} +\sum_{k=1}^{n}x_k\underbrace{E_k+E_{n+k}}_{A_k}\\
	&=&\begin{bmatrix}
		\Re(A_0)&-\Im(A_0)\\
		\Im(A_0)&\Re(A_0)
	\end{bmatrix} +\sum_{k=1}^{n}x_kA_k,
\end{eqnarray*}
However, $A(x)$ (with the same $x$) is a minimal matrix of $n\times n$. 

Next, we obtain another characterization of the solution of problem \eqref{P1} without convert it into a real SDP problem. 

\begin{theorem}\label{prop varias equiv de matriz minimal}
	Let $A\in  M_n^h(\C)$ and $x\in \R^n$. The following conditions are equivalent:
	\begin{enumerate}
	\item[a) ] $x= \Diag(A)$ is a solution of \eqref{P2} (i.e: $A=A(x)$ is minimal).
	\item[b) ] (Adapted to the more general case of $A\in M_n^h(\C)$ from \cite[Theorem 3.2]{overton}) 
	\\
	If $\{q_1,\dots, q_n\}$ is an orthonormal basis of eigenvectors of $A$ corresponding to the eigenvalues $\|A\|=\lambda_1\geq \dots\geq\lambda_n=-\|A\|$, 
	 $Q_1=[q_1|\dots|q_t]$, $Q_2=[q_{n-s+1}|\dots|q_n]$ are the $n\times t$ and $n\times s$ matrices whose columns correspond to the eigenvectors of $\lambda_1$ and $\lambda_n$ respectively, then there exist semidefinite positive self-adjoint matrices $U\in \C^{t\times t}$ and $V\in \C^{s\times s}$ such that
\begin{equation}
\label{b)(i)tr U+tr V=1} 
\tr(U)+\tr(V)=1, \text{ and }
\end{equation}			
\begin{equation}
\label{b)(ii)}  \tr(Q_1^*E_kQ_1U)-\tr(Q_2^*E_kQ_2V)=0,\ \forall k=1,...,n.
\end{equation}			
	\item[c) ] (Adapted from \cite[Theorem 2.1.6]{zhang-jiang} to the particular case of $W(H)=\Diag(M_n^h)$) 
	\\
	There exists $X\in M_n(\C)$ with  $\Diag(X)=0$ such that $A X = \|A\| |X|$, where $|X|=(X^*X)^{1/2}$.
	\item[d) ] 
	(From \cite[Theorem 2 (ii)]{alrv}) 
	Let $E_+$ (respectively $E_-$) be the spectral projection of $A$ corresponding to the eigenvalue $\lambda_1=\lambda_{max}(A)$ (respectively $\lambda_n=\lambda_{min}(A)$). Then there is a non-zero $X\in M_n^h(\C)$ such that
	$$
	\Diag(X)=0,\ E_+X^+=X^+,\ E_-X^-=X^-\ \ and \ \ \tr(A X)=\|A\|\, \|X\|_1,
	$$
	where $
	X^+=\frac{|X|+X}2\ \text{ and } \ 
	X^-=\frac{|X|-X}2 \text{ (with } |X|=(X^2)^{1/2}\geq 0\text{).}
	$
\end{enumerate}
\end{theorem}

\begin{proof} The equivalences a) $\Leftrightarrow$ c) $\Leftrightarrow$ d) have already been proved in the provided citations. The equivalence with item b) is the only that needs a proof.
\\	
	Let $W$ be the unitary $n\times n$ matrix whose columns are the eigenvectors of $A$:
	$$
	W=[Q_1|Q_2|R],
	$$
	  with  $R=[q_{t+1}|\dots |q_{n-s}]$ (following the notation used in \eqref{def Q1 y Q2}) and $Q_1$, $Q_2$ from b).
\\	
	Now consider the diagonal blocks of $t\times t$, $s\times s$ and $(n-t-s)\times(n-t-s)$ to define the following $n\times n$ self-adjoint matrix
	$$
	X=W\cdot \left(\begin{smallmatrix}U&0&0\\0&-V&0\\0&0&0 \end{smallmatrix}\right)\cdot W^*
	$$
	using the positive semidefinite matrices $U$ and $V$ from b) of sizes  $t\times t$ and $s\times s$ respectively.
\\	
Then using \eqref{b)(ii)} it can be proved that $\Diag(X)=0$. Moreover, using \eqref{b)(i)tr U+tr V=1} and the fact that $U\geq0$, $V\geq 0$ we obtain
\begin{equation*}
\begin{split}
\|X\|_1&=\tr|X|=\tr \left|W\cdot \left(\begin{smallmatrix}U&0&0\\0&-V&0\\0&0&0 \end{smallmatrix}\right)\cdot W^*\right|
= \tr\left( W \left|\left(\begin{smallmatrix}U&0&0\\0&-V&0\\0&0&0 \end{smallmatrix}\right)\right| W^*\right)
=\tr  \left|\left(\begin{smallmatrix}U&0&0\\0&-V&0\\0&0&0 \end{smallmatrix}\right)\right| \\
&=\tr  \left(\begin{smallmatrix}U&0&0\\0&V&0\\0&0&0 \end{smallmatrix}\right)=1.
\end{split}
\end{equation*}

And $A$ can be diagonalized using $W$ as
$$
A=W\cdot\left(\begin{smallmatrix}\lambda_1 I_t&0&0\\0&-\lambda_1 I_s&0\\0&0&D_0 \end{smallmatrix}\right)\cdot W^*
$$
where $D_0$ is the diagonal matrix with the eigenvalues of $A$ distinct from $\lambda_1$ and $-\lambda_1=\lambda_n$ (including multiplicity) and $I_k$ denotes the $k\times k$ identity matrix.
\\
Hence
\begin{equation*}
\begin{split}
A\cdot X&=W\cdot\left(\begin{smallmatrix}\lambda_1 I_S&0&0\\0&-\lambda_1 I_t&0\\0&0&D_0 \end{smallmatrix}\right)\cdot W^*
\cdot W\cdot \left(\begin{smallmatrix}U&0&0\\0&-V&0\\0&0&0 \end{smallmatrix}\right)\cdot W^* 
=W\cdot \left(\begin{smallmatrix}\lambda_1 U&0&0\\0&\lambda_1 V&0\\0&0&0 \end{smallmatrix}\right)\cdot W^*
\\
&=\lambda_1\ W\cdot \left(\begin{smallmatrix} U&0&0\\0& V&0\\0&0&0 \end{smallmatrix}\right)\cdot W^*=\lambda_1 |X|
\end{split}
\end{equation*}
hold, which together with the fact that $\lambda_1=\|A\|$ and item c) of Proposition \ref{prop varias equiv de matriz minimal} imply that $A$ is a minimal matrix.

To prove de implication a) $\Rightarrow$ b) we will use that the statement d) is equivalent to the condition of being a minimal matrix. Then given a minimal matrix $A\in M_n^h(\C)$ there exists $X$ such that $\Diag(X)=0$, $\tr(|X|)=\|X\|=1$, $\tr(A X)=\|A\|$, $E_+ X =X$ and  $E_- X =X$. 
\\
Now consider the unitary matrix $Q=[Q_1|Q_2|Q_3]$  constructed as follows. The  $n\times t$ and $n\times s$ matrices $Q_1=[v_1|\dots|v_t]$, $Q_2=[v_{n-s+1}|\dots|v_{n}]$ are constructed with columns of eigenvectors corresponding to the eigenvalues $\lambda_1=\|A\|$ and $\lambda_n=-\|A\|$, and $Q_3=[v_{t+1}|\dots|v_{n-s}]$ is a $n\times (n-s-t)$ matrix with columns formed by eigenvectors of $A$ that complete an orthonormal basis of $\C^n$. 
\\
Then, from the proof of (i) $\Rightarrow$ (ii) in \cite[Theorem 2 (ii)]{alrv} follows that there exists $Y\geq 0$, $Z\geq 0$ such that
$$
X=Q \left(\begin{smallmatrix}Y&0&0\\0&-Z&0\\0&0&0 \end{smallmatrix}\right)Q^*
\ \text{ and } \ 
A=Q\cdot \left(\begin{smallmatrix}\lambda_1 I_t&0&0\\0&-\lambda_1 I_s&0\\0&0&D_0 \end{smallmatrix}\right)\cdot Q^*. 
$$
Then since $\tr|X|=1$ it must be 
\begin{equation}
\label{eq tr de Y+Z=1}
1=\tr|X|=\tr \left(Q \left|\left(\begin{smallmatrix}Y&0&0\\0&-Z&0\\0&0&0 \end{smallmatrix}\right)\right|Q^*\right)=\tr   \left(\begin{smallmatrix}|Y|&0&0\\0&|Z|&0\\0&0&0 \end{smallmatrix}\right) =\tr(Y+Z).
\end{equation}
Moreover, using the expression  $Q=[Q_1|Q_2|Q_3]=  \begin{pmatrix}
Q_1 & Q_2 & Q_3
\end{pmatrix} $ we can compute
$$
X= \begin{pmatrix}
Q_1 & Q_2 & Q_3
\end{pmatrix} \left(\begin{smallmatrix}Y&0&0\\0&-Z&0\\0&0&0 \end{smallmatrix}\right) \left( \begin{smallmatrix}
Q_1 \\ Q_2 \\ Q_3
\end{smallmatrix}\right)^*= Q_1YQ_1^*-Q_2ZQ_2^*.
$$
Hence, since $\Diag(X)=0$ follows that 
$$
X_{k,k}=(Q_1 Y Q_1^*-Q_2 Z Q_2^*)_{k,k}=\tr(E_k (Q_1 Y Q_1^*-Q_2 Z Q_2^*) E_k)=0,\ \forall k=1,\dots, n
$$
and therefore 
\begin{equation}
\label{ec tr de la suma es 0}
0=\tr\left(E_k (Q_1 Y Q_1^* - Q_2 Z Q_2^*)\right)=\tr\left(Q_1^* E_k Q_1 Y -Q_2^*E_k Q_2 Z \right)
\end{equation}
Then considering $U=Y$ and $Z=V$ the equations \eqref{eq tr de Y+Z=1} and \eqref{ec tr de la suma es 0} prove that item b) holds if $A$ is minimal.
\end{proof}
Observe that the proof of Theorem \ref{prop varias equiv de matriz minimal} is different than the one made by Overton \cite{overton} and Fletcher \cite{fletcher}, and it is related with the characterization of minimal self-adjoint matrices made in \cite{alrv}.

\section{Subdifferential and moment of a subspace} \label{sec subdif}
Here we generalize some of the results developed by Overton in \cite{overton2} to complex matrices and compact self-adjoint operators. Also, we relate the concept of subdifferential of the maximum eigenvalue with the moment of a subspace. In order to do this, it is necessary to state some particular definitions and previous results.

Recall that, given a subspace $S$ of a separable Hilbert space $H$, the moment of $S$ is defined by
\begin{equation}\label{defin moment}
m_S=\text{co}\{|v|^2: v\in S, \|v\|=1\},
\end{equation}
where $v=(v_1,v_2,...)$ are the coordinates of $v$ any fixed basis of $H$, and $|v|^2=\left( |v_1|^2, |v_2|^2,...\right) $. In particular, if $H=\C^n$, by Proposition 3.2 in \cite{KV momento}
\begin{equation} \label{defi moment Y matrix}
m_S=\Diag\left(\{Y\in M_n^h(\C): Y\geq 0, \tr(Y)=1, Im(Y)\subset S \} \right), 
\end{equation}
and if $H$ is an infinite dimensional Hilbert space, by Proposition 1 in \cite{BV-Moment}
\begin{equation} \label{defi moment Y compact}
m_S=\Diag\left(\{Y\in B_1(H)^{sa}: Y\geq 0, \tr(Y)=1, Im(Y)\subset S \} \right). 
\end{equation}
On the other hand, consider a sequence $\mathbf{A}=\{A_j\}_{j=1}^\infty$ of self-adjoint compact operators or matrices $A_j$ with bounded norm ($\|A_j\|\leq c$, for all $j$).	We define the joint numerical range of $\mathbf{A}$ by
\begin{equation}
	\label{eq def JNR mat}
	W\left(\mathbf{A}\right)=\left\{\{\tr \left(\rho A_j\right)\}_{j=1}^\infty: \rho\in 
	M_n(\C)^h	\wedge \tr(\rho)=1 \wedge \rho\geq 0\right\},
\end{equation}
when $\mathbf{A}\subset M_n(\C)^h$, and
\begin{equation}
		\label{eq def JNR}
		W\left(\mathbf{A}\right)=\left\{\{\tr \left(\rho A_j\right)\}_{j=1}^\infty: \rho\in 
		\mathcal{B}_1^{sa}(H)	\wedge \tr(\rho)=1 \wedge \rho\geq 0\right\},
\end{equation}
when $\mathbf{A}\subset \mathcal{K}^{sa}(H)$.

For $E=\{e_j\}_{j=1}^{\infty}$ we will denote with  $e_j\otimes e_j=E_j$, the rank-one orthogonal projections onto the subspaces generated by $e_j\in E$, for all $j\in\mathbb{N}$.
We will be particularly interested in the study of $W(\mathbf{A})$ in the case of 
$\mathbf{A}=\mathbf{A_{S,E}}=\{  P_SE_j P_S \}_{j=1}^\infty$ and $S$ a finite dimensional subspace of $H$ 
\begin{equation}\label{def W AsubS,E}
	W\left(\mathbf{A_{S,E}}\right)=\left\{\{\tr \left(P_SE_j P_S\rho\right)\}_{j=1}^\infty: \rho\in \mathcal{B}_1(H) ,\ \rho\geq 0 \text{ and } \tr(\rho)=1\right\}.
\end{equation}
Note that $W\left(\mathbf{A_{S,E}}\right)\subset \ell^1\left(\R\right)\cap \mathbb{R}_{\geq 0}^\mathbb{N}$.
In this context, we will consider the set of its density operators 
\begin{equation}
	\label{def DsubS}
	\mathcal{D}_S=
	\left\{Y\in \mathcal{B}_1(H):P_S Y=Y\geq 0\ ,\text{tr}(Y)=1\right\}
\end{equation}
(note that $P_S Y=YP_S=P_S YP_S$ for $Y\in\mathcal{D}_S$). If $\dim S<\infty$, the affine hull of $\mathcal{D}_S$ is also finite dimensional.

There exists a relation between the moment of a subspace and the Joint numerical range of the particular family $\mathbf{A_{S,E}}$, as we illustrate in the next result. 
\begin{proposition}[Proposition 2 \cite{BV-Moment}]
	\label{prop: equivalencias de momento}
	The following are equivalent definitions of $m_S$, the moment of $S$ (see \eqref{def W AsubS,E}), with $\dim S=r <\infty$, related to a basis $E=\{e_i\}_{i=1}^\infty$ of $H$. Note the identification made between diagonal operators and sequences.
	\begin{enumerate}
		\item[a) ] \label{defequiv mS 1 en prop} $m_S=\text{Diag}(\mathcal{D}_S)$.
		\item[b) ] \label{defequiv mS 2 en prop}
		$
		m_S= \text{co}\left\{|v|^2: v \in S \text{ and } \|v\|=1 \right\}.
		$
		\item[c) ] \label{defequiv mS 3 en prop} $m_S= \bigcup\limits_{\{s^i \}_{i=1}^r  \text{o.n. set in } S} \  \text{co} \{|s^i|^2 \}_{i=1}^r.
		$
		\item[d) ] \label{eq relac momento y trEiY} 
		$m_S=\{\left(\tr(E_1 Y),\dots,\tr(E_n Y),\dots\right)\in\ell^1(\R): Y \in \mathcal{D}_S \}$.
		\item[e) ] \label{eq relac momento y JNR} 
		$m_S=W(P_S E_1 P_S,\dots,P_S E_n P_S,\dots)\cap \left\{x\in \ell^1(\R):x_i\geq 0 \text{ and } \sum_{i=1}^\infty x_i=1\right\}$, 
		where $P_S$ is the orthogonal projection onto $S$, and $W$ is the joint numerical range 
		\eqref{eq def JNR}.
	\end{enumerate}
\end{proposition}

See \cite{KV momento, BV-Moment} for more properties about $m_S$ in finite and infinite dimensional cases, respectively.

\subsection{The finite dimensional case}

Let $A\in M_n^h(\C)$, define $\lambda_{max}(A)=\lambda_1(A)$ and assume that it has multiplicity $s\geq 1$. Then, $\lambda_1: M_n^h(\C)\to \R $ is a convex function, since it can be written as the maximum of a set of linear functions,
\begin{eqnarray} \label{ray mat}
	\lambda_1(A)&=&\max\{\left\langle Aq,q\right\rangle:q\in \C^n, \|q\|=1\}=\max\ \{\tr(Aqq^*): q\in \C^n, \|q\|=1\}\\ \nonumber
	&=&\max \{\left\langle A,R \right\rangle_{tr}:\ R\in M_n^h(\C), U\geq 0, \tr(R)=1\}. 
\end{eqnarray}
The proof of the previous fact is done in Proposition \ref{lemma rayleigth} in a more general context. 

\begin{definition}\label{def subdiferencial}
	For any convex function $f:\mathcal{X}\to \R$ defined on a Banach space $\mathcal{X}$, and $\mathcal{X}^*$ its dual, it can be defined the \textbf{subdifferential} of $f$ at $x\in\mathcal{X}$ as
\begin{equation}\label{subdif gral banach}
\partial f(x)=\{v\in\mathcal{X}^*: f(y)-f(x)\geq \text{Re } v(y-x), \forall\ y\in \mathcal{X}\}, 
\end{equation}
as in \cite{bata-grover}. 
\end{definition}
In particular, if $\mathcal{X}=M_n^h(\C)$ and $f(\cdot)=\lambda_1(\cdot)$, the subdifferential at $x=A\in M_n(\C)^h$ is
$$\partial \lambda_1(A)=\{V\in M_n^h(\C): \lambda_1(Y)-\lambda_1(A)\geq \text{Re } \left\langle V,(Y-A)\right\rangle_{tr},\ \forall\ Y\in M_n^h(\C) \}.
$$
Then, using \eqref{ray mat} and  similar arguments than those in \cite{overton2}, the subdifferential of $\lambda_1$ at $A$ is the set
\begin{eqnarray} \label{subdif matrices}
\partial\lambda_1(A)&=&\text{co}\{qq^*: Aq=\lambda_1(A)q \text{ and }\|q\|=1\}\\ \nonumber
&=&\{Q_1R_sQ_1^*: R_s\in M_s^h(\C), R_s\geq 0, \tr(R_s)=1\},
\end{eqnarray}
where the columns of $Q_1$ form an orthonormal set of $s$ eigenvectors for $\lambda_1(A)$ (an orthonormal basis of the eigenspace of $\lambda_1(A)$). Observe that $Q_1$ depends on the matrix $A$.

Let $A(x)=A_0+\Diag(x)=A_0+\sum_{j=1}^nx_ke_ke_k^*$, with $A_0\in  M_n^h(\C)$, $\{e_k\}_{k=1}^{n}$ a fixed orthonormal basis of $\C^n$ and $x\in \R^n$. The maximum eigenvalue of $A(x)$, $\lambda_1(A(x))=\lambda_{max}(A(x))$, is a map from $\R^n$ to $\R$. Observe that $\lambda_1(A(x))$ is a composition of a smooth function $A(\cdot)$ and a convex map $\lambda_1(\cdot)$. Moreover, for every $k$, the partial derivatives of $A$ are
$$\frac{\partial A}{\partial x_k}(x)=e_ke_k^*.$$
Adapting Theorem 3 in \cite{overton2} to the self-adjoint case, the subdifferential of $\lambda_1(A(x))$ is
\begin{equation} \label{subdif A(x)}
	\partial\left(\lambda_1(A(x))\right)=\{v\in \R^n:\ v_k=\left\langle R_s, Q_1(x)^*e_ke_k^*Q_1(x)\right\rangle_{tr}, R_s\in M_s^h(\C), R_s\geq 0, \tr(R_s)=1\},
\end{equation}
where $s$ is the multiplicity of $\lambda_1(A(x))$ and the columns of $Q_1(x)$ form an orthonormal basis of eigenvectors for $\lambda_1(A(x))$. Note that $Q_1(x)$ depends on $A(x)$ and
$$
v_k=\left\langle R_s, Q_1(x)^*e_ke_k^*Q_1(x)\right\rangle_{tr}=\tr(R_s Q_1(x)^*e_ke_k^*Q_1(x))\geq 0,$$
for every $k$, since $R_s$ and $Q_1(x)^*e_ke_k^*Q_1(x)$ are semidefinite positive matrices. 

Using \eqref{subdif A(x)}, we obtain the following characterization of the subdifferential.

\begin{theorem} \label{teo subdif=moment}
Let $A(x)=A_0+\Diag(x)$, with $A_0\in  M_n^h(\C)$ and $x\in \R^n$, and $S_1$ be the eigenspace of $\lambda_1(A(x))$. Then
\begin{equation}\label{subdiff l1 = momento}
 \partial\big(\lambda_1(A(x))\big)=\Diag(\partial\lambda_1(A(x)))= m_{S_1},
\end{equation}
where $m_{S_1}$ is the moment of the eigenspace $S_1$ and we have identified diagonal matrices with vectors in the last equality.
\end{theorem}

\begin{proof}
Suppose $\lambda_1(A(x))=\lambda_{max}(A(x))$ has multiplicity $s$ and $S_1$ is the eigenspace of $\lambda_1(A(x))$ with a fixed orthonormal basis of eigenvectors $\{q_1(x),...,q_s(x)\}$. If $Q_1(x)=\begin{bmatrix}
	q_1(x)|...|q_s(x)
\end{bmatrix}\in \C^{n\times s}$.  By \eqref{subdif A(x)}, any $v\in\partial\lambda_1(A(x))$ has coordinates
$$v_k=\tr \left( R_s Q_1(x)^*e_ke_k^*Q_1(x)\right)=\tr \left(Y(x)e_ke_k^*\right)=\left\langle Y(x)e_k,e_k\right\rangle=Y_{kk},$$ 
with $R_s\in M_s^h(\C)$, $R_s\geq 0$, $\tr(R_s)=1$, for every $k=1,...,n$. Then, $v=(v_1,v_2,...,v_n)=\Diag(Y(x))$ with $Y(x)=Q_1(x)R_s Q_1(x)^*\in \partial \lambda_1(A(x))$ satisfies
\begin{itemize}
	\item $Y(x)=Y(x)^*\geq 0$,
%
	\item $\tr(Y(x))=\tr(Q_1(x)R_s Q_1(x)^*)=\tr(R_sQ_1(x)^*Q_1(x))=\tr(R_s)=1$, and
	
	\item $\text{Im}(Y(x))\subset S_1$. Indeed, for every $h\in \C^n$ note that
	$$Y(x)h=Q_1(x)\left( \underbrace{R_s Q_1(x)^*h}_w \right)=Q_1(x)w\in S_1,$$
	with $w$ a column vector of $s$ coordinates. 
	
\end{itemize}
Therefore, by \eqref{defi moment Y matrix}
$$v=\Diag(Y(x))=\Diag\left( Q_1(x)R_s Q_1(x)^*\right) \in m_{S_1}.$$

On the other hand, take any $v\in m_{S_1}$. Then, it can be written as
$$
v=\left( \tr(e_1e_1^*Y), \tr(e_2e_2^*Y), ...,\tr(e_ne_n^*Y) \right), $$
with $Y=Y^*\geq 0$, $\tr(Y)=1$ and $Im(Y)\subset S_1$.
In terms of the orthogonal decomposition $\C^n=S_1\oplus S_1^{\perp}$, given by the matrix $Q=\begin{bmatrix}
Q_1(x)&Q_2(x)
\end{bmatrix}$ ($Q_2$ is a matrix whose columns form an orthonormal set for $S_1^{\perp}$ and $Q$ is an unitary matrix), $Y$ is defined by
$$Y=Q\begin{bmatrix}
V&0\\
0&0
\end{bmatrix}Q^*=Q_1(x)VQ_1(x)^*$$
with $V\in M_s^h(\C)$, $1=\tr(Y)=\tr(V)$ and $V\geq 0$. Therefore, $v\in \partial \lambda_1(A(x))$.
\end{proof}

\begin{corollary}
Under the assumptions of Theorem \ref{teo subdif=moment}, if $\lambda_1(A(x))$ has multiplicity one (i.e., $s=1$), then 
$$\partial \lambda_1(A(x))=\{|v|^2:A(x)v=\lambda_1(A(x))v,\|v\|=1\},$$
$\lambda_1(A(x))$ is derivable and $\frac{\partial \lambda_1}{\partial x_k}(x)=|v_k|^2$ for every $k=1,2,\dots n.$
\end{corollary}

The next result gives a concrete formula to the directional derivative of $\lambda_1(A(x))$ and appears in \cite{overton2}, but here we include an explicit proof for the self-adjoint case. 

\begin{proposition}
Let $A(x)=A_0+\Diag(x)$, with $A_0\in  M_n^h(\C)$ and $x\in \R^n$.	Suppose $\lambda_1(A(x))=  \lambda_{max}(A(x))$, has multiplicity $s$, with
a corresponding orthonormal basis of eigenvectors $\{q_1(x),...,q_s(x)\}$ and $Q_1(x)=\begin{bmatrix}
	q_1(x)|...|q_s(x)
\end{bmatrix}$. Then the directional derivative of $\lambda_1$ at $x\in \R^n$ in the direction $w\in \R^n$ that is defined by
$$\lambda'_1(x,w)=\lim\limits_{t\to 0^+}\dfrac{\lambda_1(x+tw)-\lambda_1(A(x))}{t}$$
is the largest eigenvalue of 
$$B(w)=\sum_{k=1}^{n}w_kQ_1(x)^*e_ke_k^*Q_1(x).$$
\end{proposition}

\begin{proof}
Recall that $\lambda_1(A(x))=\lambda_1\circ A(x)$, is a composition of a smooth map $A(x)$ with a convex function $\lambda_1$. Then, for every $w\in \R^n$ 
$$\lambda_1'(x,w)=\max_{v\in \partial \lambda_1(A(x))} \left\langle v,w\right\rangle,$$
since the generalized derivative and generalized gradient coincide with the directional derivative and subdifferential, respectively (see Proposition 2.2.7 in \cite{clarke1}). Therefore,
\begin{eqnarray*}
\lambda_1'(x,w)&=&\max\left\lbrace \sum_{k=1}^{n} v_kw_k: v_k=\left\langle R_s, Q_1(x)^*e_ke_k^*Q_1(x)\right\rangle_{tr},R_s\in M_s^h(\C), R_s\geq 0, \tr(R_s)=1\right\rbrace\\
&=&\max\left\lbrace  \left\langle R_s,\sum_{k=1}^{n}w_k Q_1(x)^*e_ke_k^*Q_1(x)\right\rangle_{tr}:\ R_s\in M_s^h(\C), R_s\geq 0, \tr(R_s)=1\right\rbrace\\
&=&\max\left\lbrace  \left\langle R_s,B(w)\right\rangle_{tr}:\ R_s\in M_s^h(\C), R_s\geq 0, \tr(R_s)=1\right\rbrace\\
&=&\lambda_1(B(w)),
\end{eqnarray*} 
where the last equality is due to \eqref{ray mat}.
\end{proof}

\begin{lemma} \label{lema subdif lambda n}
Let $A(x)=A_0+\Diag(x)$, with $A_0\in M_n^h(\C)$ and $x\in \R^n$, $\lambda_n(A(x))$ be the minimum eigenvalue of $A(x)$ and $S_n$ its corresponding eigenspace. Then, 
\begin{equation}\label{subdif lambda_n}
\partial \lambda_n(A(x))=\partial \lambda_n(x)=-m_{S_n}.
\end{equation}
\end{lemma}

\begin{proof}
Since $\lambda_n(A(x))=-\lambda_1(-A(x))$ for any $A(x)\in M_n^h(\C)$, then
\begin{equation*}
\begin{split}
\partial \lambda_n(A(x))&=-\partial \lambda_1(-A(x))\\
&=-\Diag\left(\text{co} \{uu^*: \|u\|=1,\ -A(x)u=\lambda_1(-A(x))u\} \right) \\
&=-\Diag\left(\text{co} \{uu^*: \|u\|=1,\ A(x)u=\lambda_n(A(x))u\} \right) \\
&=-m_{S_n}.
\end{split}
\end{equation*}
\end{proof}

The subdifferential of the spectral norm of a matrix $A$ is 
\begin{equation}\label{watson}
	\begin{split}
	\partial\|A\|&=\text{co}\{uv^*: Au=\|A\|v \text{ and }\|v\|=1\}\\
	&=\{VR_{t}W^*:R_{t}\in M_{t}^h(\C), R_{t}\geq 0, \tr(R_{t})=1\},
	\end{split}
\end{equation}
where $V, W$ are unitary matrices of the singular value decomposition of $A$, $A=V\D(s(A))W^*$ and 
$s_1(A)=\|A\|$. The proof of \eqref{watson} appears first in \cite{watson} for the real case and, more recently for the complex case, in \cite{bata-grover} and \cite{Grover}.

The subdifferential \eqref{watson} can be closely related with the subdifferentials of $\lambda_1$ and $\lambda_n$ in some cases, as we observe in the next statement.

\begin{remark} \label{eq subdif norma -eig matrix}
The expression in \eqref{watson} for any $A(x)=A_0+\Diag(x)$, with $A\in M_n^h(\C)$ and $x\in \R^n$ is
\begin{equation}\label{watsonhermit}
\partial\|A(x)\|=\text{co}\{uu^*: A(x)u=\|A(x)\|u \text{ and }\|u\|=1\}.
\end{equation}
Let $\lambda_n(A(x))$ and $\lambda_1(A(x))$ be the minimum and maximum eigenvalue of $A(x)$, respectively. Considering \eqref{watsonhermit} and Lemma \ref{lema subdif lambda n}, it is evident that 
$$\partial\|A(x)\|=\left\lbrace \begin{array}{lll}
\partial\lambda_1(A(x))&\text{if}& \|A(x)\|=\lambda_1(A(x))\\
\partial\lambda_n(A(x))&\text{if}& \|A(x)\|=-\lambda_n(A(x))=|\lambda_n(A(x))|\\
{\rm co}\left(\partial \lambda_1(A(x))\cup\partial \lambda_n(A(x)) \right)&\text{if}& \|A(x)\|=-\lambda_n(A(x))=\lambda_1(A(x)).
\end{array}\right. $$

\end{remark}

%
%
%
%
%
%
%

\begin{theorem}\label{teo equiv subdif, momentos, etc para matrices}
Let $A(x)=A_0+\Diag(x)$, with $A_0\in  M_n^h(\C)$ and $x\in \R^n$ such that $\lambda_1(A(x))=-\lambda_n(A(x))$. Then, the following statements are equivalent,

\begin{enumerate}
	\item[(1)] $0\in \partial\|A(x)\|$.
	\item[(2)] $0\in\partial\lambda_1(A(x))+\partial\lambda_n(A(x))$.
	\item[(3)] $m_{S_1}\cap m_{S_n}\neq \emptyset$, where $S_1$ and $S_n$ are the eigenspaces of $\lambda_1(A(x))$ and $\lambda_n(A(x))$, respectively. 
	\item[(4)] $W\left( \left\lbrace P_{S_1}e_ie_i^*P_{S_1}\right\rbrace_{i=1}^n \right)\cap W\left( \left\lbrace P_{S_n}e_ie_i^*P_{S_n}\right\rbrace_{i=1}^n \right)\neq \{0\} $. 
	\item[(5)] $A(x)$ is minimal.
\end{enumerate}
\end{theorem}

\begin{proof}
	
	
	\vspace{1cm}
	
The equivalences $(3) \Leftrightarrow (4) \Leftrightarrow (5)$ have already been proved in \cite{KV momento}. 

(1)$\Leftrightarrow$(3) If $0\in \partial \|A(x)\|$ and $\lambda_1(A(x))=-\lambda_n(A(x))$, then using Remark \ref{eq subdif norma -eig matrix}
$$0\in {\rm co}\left(\partial \lambda_1(A(x))\cup\partial \lambda_n(A(x)) \right)={\rm co}\left(m_{S_1}\cup -m_{S_n} \right)
$$ 
and there exist $\alpha\in (0,1)$, $Y_0\in \{Y\in M_n^h(\C): Y\geq 0, \tr(Y)=1, \text{Im}(Y)\subset S_1\}$ and $Z_0\in \{Z\in M_n^h(\C): Z\geq 0, \tr(Z)=1, \text{Im}(Z)\subset S_n\}$ such that $0=\alpha \D(Y)+(1-\alpha)\D(-Z)$. Using that $\tr(Y)=\tr(Z)=1$ we obtain that $\alpha=\frac 12$ and then $\D(Y)=\D(Z)$. Therefore, $m_{S_1}\cap m_{S_n}\neq \emptyset$. 

The converse implication can be proved reversing the previous steps.

%
%

To prove (2)$\Leftrightarrow$(3) we can use 
the formulas $\partial\big(\lambda_1(A(x))\big)= m_{S_1}$
and $\partial \lambda_n(A(x))=-m_{S_n}$ from  \eqref{subdiff l1 = momento} and  \eqref{subdif lambda_n}. Then it is trivial that   $m_{S_1}\cap m_{S_n}\neq \emptyset$ if and only if $0\in m_{S_1}-m_{S_n}=\partial \lambda_1(A(x))+\partial \lambda_n(A(x))$.
\end{proof}

\subsection{The compact operator case}
\begin{lemma}\label{rango 1 y tipo traza}
Let $B_1(H)$ be the ideal of trace class operators. Then,
$${\rm co}\left(\{hh^*: h\in H, \|h\|=1\} \right)=\{Y\in B_1(H): Y\geq 0,\ \tr(Y)=1 \},$$
\end{lemma}
\begin{proof}
If $\sum_{j}a_jh^j(h^j)^*$ is a convex combination of unitary vectors $h^j\in H$, then it fulfills that is a semidefinite positive compact operator with $\tr\left( \sum_{j}a_jh^j(h^j)^*\right) =\sum_{j}a_j=1$. On the other hand, every $Y\in B_1(H)^{sa}$, with $Y\geq 0$ and $\tr(Y)=1 $ can be written as a (maybe infinite) convex combination of rank one operators .

\end{proof}

\begin{definition}\label{defi reg}
	Given a Banach space $\mathcal{X}$, a function $f: \mathcal{X}\to \R$ is said to be regular at $x\in \mathcal{X}$ if
	\begin{enumerate}
		\item for all $v$, the usual one-sided derivative 
		$$f'(x,v)=\lim\limits_{t\to 0^+}\dfrac{f(x+tv)-f(x)}{t}$$
		exists.
		\item For all $v$, $f'(x,v)$ coincides with the general derivative.
	\end{enumerate}
\end{definition}
To see more details of this definition, see \cite{clarke2} and  \cite{clarke1}.

\begin{proposition} \label{lemma rayleigth}
Let $A\in K(H)^{sa}$ and 
\begin{equation} \label{eig function}
 \lambda_{max}(A)=\max\{\lambda\in \C:\ A-\lambda I \text{ is not invertible}\}
\end{equation}
be the maximum eigenvalue of the spectrum of $A$, then the following statements hold.
\begin{enumerate}
	\item $ \lambda_{max}(A)\in \R$ is an eigenvalue of $A$ and it has finite multiplicity if $ \lambda_{max}(A)\neq 0$.
	\item The following are equivalent forms to describe $ \lambda_{max}(A)$, 
	\begin{eqnarray}
	\label{ray1}
	\lambda_{\max}(A)&=&\max_{\|h\|=1}\left\langle Ah,h \right\rangle \\  \label{ray2}
	&=&\max\{\tr(Ahh^*): \|h\|=1\}=\max\{\left\langle A,hh^*\right\rangle_{tr}:\ \|h\|=1\}\\
	&=&\max \{\left\langle A,Y \right\rangle_{tr}:\ Y\in B_1(H)^{sa},Y\geq 0, \tr(Y)=1\}. \label{rayleight}
	\end{eqnarray}
	\item $\lambda_{\max}:K(H)^{sa}\to \R$ is a convex function and is Lipschitz near $A$ and regular in the sense of Definition \ref{defi reg}.
	\item As a particular case of Lemma \ref{rango 1 y tipo traza},
	we define the set
	\begin{equation} \label{defi d0}
	\mathcal{D}_ {S_{\max}}= \text{co}\left(\{qq^*: q\in H, Aq\right. = \left.\lambda_{\max}(A)q, \|q\|=1\} \right).
	\end{equation}
	Then, 
	\begin{equation}\label{convex hull and trace class for eigespace}
	\mathcal{D}_ {S_{\max}}= \{Y\in B_1(H)^{sa}: Y\geq 0 , E_+Y=YE_+,\ \tr(Y)=1 \},
	\end{equation}		
	where $E_+$ is the orthogonal projection onto the eigenspace of $\lambda_{\max}(A)$. Moreover, if $\lambda_{\max}(A)\neq 0$, then
	 \begin{equation} \label{convex hull and trace class finite mult}
	 \mathcal{D}_ {S_{\max}}=\{ R_sQ_{\max}^*: R_s\in M_s^h(\C),  R_s\geq 0, \tr(R_s)=1\},
	 \end{equation}
 where $s$ is the multiplicity of $\lambda_{\max}$ and the columns of $Q_{\max}$ form an orthonormal basis of eigenvectors of $\lambda_{\max}$.
\end{enumerate} 

\end{proposition}

\begin{proof}
\begin{enumerate}
	\item  It is a well-known fact of the spectrum of compact operators.
	\item Equality \eqref{ray1} holds since for any $\lambda$ eigenvalue  of $A$ and $v\in S_{\lambda}$, $v\neq 0$, $Av=\lambda_{\max}(A)v$ and $\left\langle Av,v \right\rangle=\lambda\left\langle v,v \right\rangle\in \R$. Then, 
	$$\lambda_{\max}(A)\geq \max\{\lambda\in\R: \exists v\in H \text{ such that } Av=\lambda v \}= \max\left\lbrace \dfrac{\left\langle Av,v \right\rangle}{\left\langle v,v \right\rangle}: v\in H\right\rbrace=\max_{\|h\|=1}\left\langle Ah,h \right\rangle. $$
	\eqref{ray2} follows from the equality $\left\langle Ah,h \right\rangle=\tr(Ahh^*)$, and \eqref{rayleight}  is due to Lemma \ref{rango 1 y tipo traza}, since maximizing a linear function (the trace) over a set gives the same result as maximizing it over its convex hull.
	\item By \eqref{ray1}, if $A,B\in K(H)^{sa}$ and $t\in [0,1]$, then
	\begin{eqnarray*}
	\lambda_{\max}(tA+(1-t)B)&=&\max_{\|h\|=1}\left\langle (tA+(1-t)B)h,h\right\rangle =\max_{\|h\|=1}\left[ t \left\langle Ah,h\right\rangle+(1-t) \left\langle Bh,h\right\rangle \right] \\
	&\leq&t\max_{\|h\|=1}\left\langle Ah,h\right\rangle+(1-t)\max_{\|h\|=1} \left\langle Bh,h\right\rangle\\
	&=&t\lambda_{\max}(A)+(1-t)\lambda_{\max}(B).
	\end{eqnarray*}
Therefore, $\lambda_{\max}:K(H)^{sa}\to \R$ is a convex function. On the other hand, $\lambda_{\max}$ is bounded above on a neighborhood of $A=A^*$ (since $\lambda_{\max}(A)\leq \|A\|<\infty$ for all $A\in B(H)$), so by \cite{clarke1} (Prop. 2.2.6 and 2.3.6), $\lambda_{\max}$ is Lipschitz near $A$ and regular.
\item The first equality is evident, since any $Y\geq 0$, such that $E_+Y=YE_+$ and $\tr(Y)=1$ can be written as
$$Y=\sum_{i=1}^{s}a_iqq^*,$$
where $\sum_{i=1}^{s}a_i=1$, $a_i\geq 0$ for every $i$, and  $\{q_i\}_{i=1}^s$ is an orthonormal basis of the eigenspace $S_{\max}$ of $A$. If $\lambda_{\max}(A)\neq 0$, then $s<\infty$ and we can define $Q_{\max}=[q_1|q_2|...|q_s]$ and $R_s=\Diag \left( \{a_i\}_{i=1}^s\right) \in M_s^h(\C)$, such that
$$Y=Q_{\max}R_sQ_{\max}^*.$$

\end{enumerate}
\end{proof}
For $\lambda_{\max}:K(H)^{sa}\to \R$, it can be defined the subdifferential at $A\in K(H)^{sa}$, using \eqref{subdif gral banach}, as
\begin{equation}\label{subdif eig compact}
\partial \lambda_{\max}(A)=\{Y\in B_1(H)^{sa}: \lambda_{\max}(B)-\lambda_{\max}(A)\geq \text{Re } \tr (Y(B-A)), \forall\ B\in K(H)^{sa}\}, 
\end{equation}
In the next result, we obtain more useful expressions of $\partial \lambda_{\max}$.
\begin{proposition}
If $A\in K(H)^{sa}$, $\lambda_{max}(A)=\lambda_{\max}(A)$ has multiplicity $s\geq 1$ and $\mathcal{D}_{S_{\max}}$ is as in \eqref{defi d0}, then the subdifferential of $\lambda_{\max}(A)$ is the set 
\begin{eqnarray}\label{subdif eig compact2} 
	\partial\lambda_{\max}(A)&=&\mathcal{D}_{S_{\max}}\\ \nonumber
	&=&\{Y\in B_1(H)^{sa}: Y\geq 0, E_+Y=YE_+,\ \tr(Y)=1\},
\end{eqnarray}
where $E_+$ is the orthogonal projection onto the eigenspace of $\lambda_{\max}(A)$.
In particular, if $\lambda_{\max}(A)\neq 0$, then
\begin{equation} \label{subdif eig compact2 finite}
\partial\lambda_{\max}(A)=\{Q_{\max}R_sQ_{\max}^*: R_s\in  M_s^h(\C), R_s\geq 0, \tr(R_s)=1\},
\end{equation}
where the columns of $Q_{\max}$ form an orthonormal basis of eigenvectors for $\lambda_{\max}$.
\end{proposition}

\begin{proof}
	As a consequence of \eqref{ray2}, the subdifferential of $\lambda_{\max}$ at $A$ can be expressed as
	$$\partial\lambda_{\max}(A)=\text{co}\{qq^*: Aq=\lambda_{\max}(A)q \text{ and }\|q\|=1\}.$$
	Then the formulations of $\partial \lambda_{\max}(A)$ in \eqref{subdif eig compact2} and \eqref{subdif eig compact2 finite} follow directly from \eqref{convex hull and trace class for eigespace} and \eqref{convex hull and trace class finite mult}, respectively.

\end{proof}

\begin{definition} \label{stric diff}
Let $\mathcal{X}$ and $\mathcal{Y}$ be Banach spaces. A function $F: \mathcal{X}\to \mathcal{Y}$ is strictly differentiable at $x\in \mathcal{X}$ if there exists a continuous linear operator from $\mathcal{X}$ to $\mathcal{Y}$, denoted by $D_sF(x)$, such that
\begin{equation}
	\label{def DsF}
		\lim\limits_{x'\to x, t\to 0^+}\frac{F(x'+tv)-F(x')}{t}=\text{Re}\ \tr\left(  D_sF(x),v\right) ,	
\end{equation}

for every $v\in\mathcal{X}$. The operator $D_sF(x)$ is the strict differential of $F$ at $x$.
\end{definition}

\begin{lemma} \label{$A$ properties}
Let $c_0(\R)=c_0$ be the space of real sequences that converge to $0$ and 
\begin{equation} \label{defi A(x)}
A(x)=A_0+\Diag(x)=A_0+\sum_{j\in \N}x_ke_ke_k^*
\end{equation}
be an affine function with $A_0\in K(H)^{sa}$ fixed and $x\in c_0$. 

\begin{enumerate}
	\item For every $k$,
	$$\frac{\partial A}{\partial x_k}(x)=e_ke_k^*.$$
	 and $A(\cdot)$ is a smooth function
	\item $A:c_0\to K(H)^{sa}$ is strictly differentiable at $x$ and $$D_sA(x)=\sum_{j\in \N}x_ke_ke_k^*=\D\left(x \right)\in D(K(H)^{sa}),$$
	where  $D_sA$ is the map defined in \eqref{def DsF}.
	\item $D_sA: c_0\to K(H)^{sa}$ satisfies that its adjoint $D_{s}A^*: K(H)^{sa}\to c_0$, is 
	$$D_sA(C)=\Diag(C)=\Diag\left(\left\lbrace \left\langle Ce_i,e_i \right\rangle \right\rbrace_{i\in \N} \right), \text{ for every } C\in K(H)^{sa}.$$

\end{enumerate}
\end{lemma}

\begin{proof}
The proof of item $1$ is direct, since each partial derivative of $A$ is a constant function. Then, for every $x\in c_0$ the differential $D_{s}A$ is
$$D_{s}A(x)=\sum_{j\in \N}x_ke_ke_k^*=\D\left(x \right)\in D(K(H)^{sa}).$$
Additionally, if $A$ is a smooth function, then it is strictly differentiable (\cite{clarke1}, p. 32) and $D_{s}A(x)$ is the strict derivative of $A$ at $x$. The adjoint $D_{s}A^*: K(H)^{sa}\to c_0$ fulfills
$$D_{s}A^*(C)x=\text{Re}\ \tr(C^*D_{s}A(x)),\ \forall\ C\in K(H)^{sa},\ \forall\  x\in c_0.$$
Then, for each $E_{ij}=e_ie_j^*$
$$D_{s}A^*(E_{ij})x=\text{Re}\ \tr(E_{ij}D_{s}A(x))=\Diag(x),$$
and for every $C\in K(H)^{sa}$ and $e_i$
$$D_{s}A^*(C)e_i=\text{Re}\ \tr(CD_{s}A(e_i))=C_{ii}.$$
Basically, $D_{s}A^*$ is the pinching operator, which extracts the main diagonal of every $C\in K(H)^{sa}$ (respect on the orthonormal prefixed basis $\{e_i\}_{i\in \N}$ of $H$), that is 
$$D_{s}A^*(C)=\Diag(C)=\Diag\left(\left\lbrace \left\langle Ce_i,e_i \right\rangle \right\rbrace_{i\in \N} \right).$$
\end{proof}
We are now in position to state one of the the main results of this subsection.

%

\begin{theorem} \label{teo subdif y moment compact}	
Let $\lambda_{\max}:K(H)^{sa}\to \R$ and $A: c_0\to \R$ be the functions defined in \eqref{eig function} and \eqref{defi A(x)}, respectively. 
Consider the composition map $\lambda_{\max}\circ A:c_0\to \R$, given by $\lambda_{\max}\circ A(x)=\lambda_{\max}(A(x))$. Let $s$ be the multiplicity of $\lambda_{\max}(A(x))$ and $S_{\max}$ the eigenspace of $\lambda_{\max}(A(x))$.

Then, the subdifferential of $\lambda_{\max}(A(x))$ at $x\in c_0$ is
\begin{equation}
	\label{subdif compacto} 
\begin{split}
	\partial\big(\lambda_{\max}(A(x))\big)&=\Diag \left(\mathcal{D}_{S_{\max}}\right)\\
	& =\Diag(\partial \lambda_{\max}(A(x)))\\
	&=m_{S_{\max}},
	\end{split}
\end{equation}
where $m_{S_{\max}}$ is the moment of the eigenspace $S_{\max}$ (see \eqref{defin moment}).

In particular, if $\lambda_{\max}(A(x))\neq 0$
\begin{equation}
	\partial\big(\lambda_{\max}(A(x))\big)=\Diag \left( \{Q_{\max}(x)R_sQ_{\max}(x)^*: R_s\in  M_s^h(\C), R_s\geq 0, \tr(R_s)=1\}\right),
\end{equation}
where the columns of $Q_{\max}(x)$ form an orthonormal set of eigenvectors for $\lambda_{\max}(A(x))$.
\end{theorem}

\begin{proof}
Let $x\in c_0$. As it was proved in Lemma \ref{$A$ properties}, $A(x)$ is a smooth function and, particularly, strictly differentiable at $x$. Furthermore, by Proposition \ref{lemma rayleigth}, $\lambda_{\max}$ is convex, Lipschitz near $A(x)$ and regular (in the sense of Definition \ref{defi reg}). Therefore, by Theorem 2.3.10 (Chain rule) and Remark 2.3.11 in \cite{clarke1}, 
$$
\partial\big( \lambda_{\max}(A(x))\big)=\partial \left( \lambda_{\max}\circ A\right) (x)=DA^*\partial \lambda_{\max}(A(x)),$$
where $DA^*$ is the adjoint of $DA$. By Lemma \ref{$A$ properties},  $DA^*: K(H)^{sa}\to c_0$ fulfills that
$$DA^*(C)=\Diag(C)=\Diag\left(\left\lbrace \left\langle Ce_i,e_i \right\rangle \right\rbrace_{i\in \N} \right).$$
By \eqref{subdif eig compact2}, 
$$
\partial\lambda_{\max}(A(x))=\mathcal{D}_{S_{\max}}=\{Y(x)\in B_1(H)^{sa}: Y(x)\geq 0, E_+Y(x)=Y(x)E_+,\ \tr(Y(x))=1\},
$$
Combining the above,
$$\partial \big(\lambda_{\max}(A(x))\big)=\Diag\left(\partial \lambda_{\max}(A(x))\right)=\Diag\left(\mathcal{D}_{S_{\max}}\right)= m_{S_{\max}},$$
where the last equality is due to Proposition \ref{prop: equivalencias de momento}.

On the other hand, if $\lambda_{\max}(A(x))\neq 0$, by \eqref{subdif eig compact2 finite}
$$
\partial \lambda_{\max}(A(x))=\{Q_{\max}(x)R_sQ_{\max}(x)^*: R_s\in  M_s^h(\C), R_s\geq 0, \tr(R_s)=1\},
$$
where $s$ is the multiplicity of $\lambda_{\max}(A(x))$, and the columns of $Q_{\max}(x)$ form an orthonormal set of eigenvectors for $\lambda_{\max}(x)$.
In this case, we obtain the following equality 
$$
\partial \big(\lambda_{\max}(A(x))\big)=DA^*\partial \lambda_{\max}(A(x))=\Diag\{Q_{\max}R_sQ_{\max}^*: R_s\in  M_s^h(\C), R_s\geq 0, \tr(R_s)=1\}. 
$$
\end{proof}

\begin{corollary}
Under the assumptions of Theorem \ref{teo subdif y moment compact}, the following formula holds
$$
\partial \lambda_{\max}(x)=\{v\in c_0: v_k=\tr \left( R_s Q_{\max}(x)^*e_ke_k^*Q_{\max}(x)\right), \forall k\in \N  \}. 
$$
\end{corollary}

\begin{proof}
By the mentioned Theorem and its proof, any $v\in \partial \lambda_{\max}(x)$ 
$$v=\Diag(Q_{\max}(x)R_sQ_{\max}(x)^*),$$
where $R_s\in  M_s^h(\C)$, $R_s\geq 0$, $\tr(R_s)=1$ and the columns of $Q_{\max}(x)$ form an orthonormal set of eigenvectors for $\lambda_{\max}(x)$. Then, the coordinates of $v$ are
$$
v_k=(Q_{\max}(x)R_sQ_{\max}(x)^*)_{kk}=\tr \left(Q_{\max}(x)R_sQ_{\max}(x)^*e_ke_k^*\right)=\tr \left( R_s Q_{\max}(x)^*e_ke_k^*Q_{\max}(x)\right). 
$$
for every $k\in \N$.
\end{proof}

Recently, in \cite{singla}, the author gave the following explicit expression for the subdifferential of the operator norm of $A\in B(H)$ such that $\text{dist}(A,K(H))<\|A\|$, 
\begin{equation}\label{subdif norm A operator}
\partial\|A\|=\overline{\text{co}}\ \{uv^*: u, v \in H,  Au=\|A\|v \text{ and }\|u\|=\|v\|=1\},
\end{equation}
where the closure of the convex hull $\overline{\text{co}}$ is in the operator norm.

%
%
%
%
When $A(x)$ is compact self-adjoint but not semi-definite, we obtain analogous results as Lemma \ref{lema subdif lambda n} and Remark \ref{eq subdif norma -eig matrix}, since $\lambda_{\max}(A(x))$ and $\lambda_{\min}(A(x))$ are real eigenvalues of $A(x)$ with finite multiplicity. We compile these facts in the next proposition and we omit the proof, which is similar to the matricial case (see Lemma \ref{lema subdif lambda n} and Remark \ref{eq subdif norma -eig matrix}). 

\begin{proposition} \label{subdif norm and eig compact op}
Let $A(x)=A_0+\Diag(x)$, with $A_0\in K(H)^{sa}$ and $x\in c_0$, be such that $A(x)$ is such that $\lambda_{\text{min}}(A(x))<0<\lambda_{\text{max}}(A(x))$. Then the following properties hold.
\begin{enumerate}
	\item If  $\lambda_\text{min}(A(x))$ is the minimum eigenvalue of $A(x)$ and $S_\text{min}$ its corresponding eigenspace, then
	\begin{equation}\label{subdif lamba_n compact op}
		\partial\big( \lambda_\text{min}(A(x))\big)=\partial \lambda_\text{min}(x)=-m_{S_\text{min}}.
	\end{equation}
	\item The equivalent expression of equation  \eqref{subdif norm A operator} in this case is
	\begin{equation}\label{watsonhermit compact op}
		\partial\left(\|A(x)\|\right)=\text{co}\{uu^*: A(x)u=\|A(x)\|u \text{ and }\|u\|=1\}.
	\end{equation}
	\item Considering \eqref{subdif lamba_n compact op} and \eqref{watsonhermit compact op}, it is evident that 	
	$$
		\partial\big(\|A(x)\|\big) =\left\lbrace \begin{array}{lll}
		\partial\lambda_{\max}(A(x))&\text{if}& \|A(x)\|=\lambda_{\max}(A(x))
		\\
		\partial\lambda_\text{min}(A(x))&\text{if}& \|A(x)\|=-\lambda_\text{min}(A(x))=|\lambda_\text{min}(A(x))|
		\\
		{\rm co}\big(\partial \lambda_{\max}(A(x))\cup\partial \lambda_\text{min}(A(x)) \big)&\text{if}& \|A(x)\|=-\lambda_\text{min}(A(x))=\lambda_{\max}(A(x)).
	\end{array}\right. $$
\end{enumerate}

\end{proposition}

\begin{proposition}\label{prop K mas norma K por I mmal}
 	Let $K\in K(H)^{sa}$ be such that $K\leq 0$, with $\dim(ker(K))=\infty$ and $Y,Z\in B_1(H)$  that satisfy  $Y  P_{S_{\lambda_{\min}(K)}}=Y$, $Z P_{\ker(K)}=Z$, $Y\geq 0$, $Z\geq 0$, $\tr(Y)=\tr(Z)=1$ and $\Diag(Y)=\Diag(Z)$. Then $A=K+\frac{\|K\|}2 I$ is a minimal operator.
\end{proposition}
\begin{proof}
	Observe that $K\leq 0$ implies that $-\lambda_{min}(K)
=|\lambda_{min}(K)|=\|K\|$ and hence  $S_{\lambda_{min}(K)}=S_{-\|K\|}$.

Now consider the spectral projection $P_{\lambda_{min}(K)}=P_{\|K\|}$ on the eigenspace $S_{\lambda_{min}(K)}=S_{-\|K\|}$ corresponding to the eigenvalue  $\lambda_{min}(K)=\|K\|$ and the orthogonal projection onto the kernel of $K$ denoted by $P_{\ker(K)}$. Note that then 
$A=-\|K\| P_{-\|K\|}+R+\frac{\|K\|}2 I$ with $R$ orthogonal to $P_{-\|K\|}$.
Hence we can obtain the following equalities
\begin{equation}\label{traza YA}
	\begin{split}
\tr(Y A) &= \tr\left(Y\left(K+\frac{\|K\|}2 I\right)\right)=\tr\left(Y P_{-\|K\|} K+\frac{\|k\|}2 Y\right)=
\tr\left(-\|K\| Y   +\frac{\|k\|}2 Y\right)\\
&=\tr\left(   \frac{-\|k\|}2 Y\right)=-\frac{\|K\|}2,
	\end{split}
\end{equation} 
 \begin{equation}\label{traza ZA}
	\begin{split}
		\tr(Z A)& = \tr\left(Z K+\frac{\|K\|}2 Z\right) =\tr\left(Z P_{\ker(K)} K +\frac{\|k\|}2 Z\right)=
		\tr\left(0   +\frac{\|k\|}2 Z\right)
		=   \frac{\|K\|}2.
	\end{split}
\end{equation}  
Now consider   $X=\frac{Z-Y}2\in B_1(H)$ with null diagonal, and define $\psi$ in the dual of $B(H)$ as $\psi(W)=\tr(X W)$. Then $\psi$ satisfies
\begin{itemize}
	\item  $\psi(XD)=0$ for every diagonal operator $D$ since $\Diag(X)=0$,
	\item $\psi(A)=\frac12\tr\left(ZK+\frac{\|K\|}2Z\right)-\frac12\tr\left( YK+\frac{\|K\|}2Y\right)=
	\frac12\left(
		\frac{\|K\|}2 -\left(-\frac{\|K\|}2 \right) \right)=\frac{\|K\|}2=\|A\|$, where we have used \eqref{traza YA} and \eqref{traza ZA}, and
	\abovedisplayskip=-\baselineskip
	\belowdisplayskip=0pt
	\abovedisplayshortskip=-\baselineskip
	\belowdisplayshortskip=0pt
	\item\begin{align*}
	 		\|X\|_1	&=
\tr(|X|)=\tr\left(\left|\frac{Z-Y}2\right|\right)=\frac12\tr\left((Z^2-YZ-ZY+Y^2)^{1/2}\right)
\\&=
\frac12\tr\left( (Z^2+Y^2)^{1/2}\right)=
\frac12\tr(Z+Y)=1
\end{align*} 

where we have used that $Y$ and $Z$ act on orthogonal subspaces.
\end{itemize}
This proves that $\psi$ is a witness of the minimality of $A$ with respect to the diagonal operators, and hence $A$ is minimal (see Section 5 and in particular Proposition 5.1 of \cite{rieffel} and Remark 9 of \cite{bot var best}).
\end{proof}
\begin{example}
	We describe here a concrete case where Proposition \ref{prop K mas norma K por I mmal} can be applied. Given $h\in H$, $\|h\|=1$, consider the rank one and hence compact operator $K=-h h^\star\leq 0$ with $|h_j|^2\leq \frac12$ and $h_j\neq 0$ for all $j\in \N$. Then $A=K+\frac12 I$ is minimal as can also be proved using Theorem \ref{teo minimal rango1}(2) since $A$ is minimal if and only if $-A$ is.
\end{example}

Now we can prove a similar result as that obtained in Theorem \ref{teo equiv subdif, momentos, etc para matrices} for matrices.

\begin{theorem}\label{equiv subdif, momentos, etc para compactos}
	Let $A(x)=A_0+\Diag(x)$, with $A_0\in  K^{sa}(H)$ and $x\in c_0$ such that $\lambda_{\max}(A(x))=-\lambda_{\min}(A(x))$. Then, the following statements are equivalent,
	\begin{enumerate}
		
		\item[(1)] $0\in \partial\big(\|A(x)\|\big)$.
		\item[(2)] $0\in\partial\lambda_{\max}(A(x))+\partial\lambda_{\min}(A(x))$.
		\item[(3)] $m_{S_{\max}}\cap m_{S_{\min}}\neq \emptyset$, where $S_{\max}$ and $S_{\min}$ are the eigenspaces of $\lambda_{\max}(A(x))$ and $\lambda_{\min}(A(x))$, respectively. 
		\item[(4)] $W\left( \left\lbrace P_{S_{\max}}e_ie_i^*P_{S_{\max}}\right\rbrace_{i=1}^\infty \right)\cap W\left( \left\lbrace P_{S_{\min}}e_ie_i^*P_{S_{\min}}\right\rbrace_{i=1}^\infty \right)\neq \{(0,\dots,0,\dots)\} $. 
		\item[(5)] There exists $m\in\N$, concrete C$^*$-isomorphisms $U_{S_{\max}}: M_{\dim(S_{\max})}\to P_{S_{\max}}B(H)P_{S_{\max}}$ and $U_{S_{\min}}: M_{\dim(S_{\min})}\to P_{S_{\min}}B(H)P_{S_{\min}}$ (for example, as defined in Proposition 10 and subsection 5.1 of \cite{BV-Moment})), matrices 	$\{B_j\}_{j=1}^m\subset M^h_{\dim(S_{\max})}(\C)$ and $\{C_j\}_{j=1}^m\subset M^h_{\dim(S_{\min})}(\C)$ with $B_j=U_{S_{\max}}^{-1}(P_{S_{\max}}e_j e_j^*P_{S_{\max}})$, $C_j=U_{S_{\min}}^{-1}(P_{S_{\min}}e_j e_j^*P_{S_{\min}})$ such that 
		$$
		W(\{B_j\}_{j=1}^m)\cap W(\{C_j\}_{j=1}^m)\neq \{(0,\dots,0)\}.
		$$
		\item[(6)] $A(x)$ is minimal.
	\end{enumerate}
\end{theorem}
\begin{proof}
	
	The equivalences $(3) \Leftrightarrow (4) \Leftrightarrow (5) \Leftrightarrow (6)$ have already been proved in \cite[Proposition 12 and Theorem 4]{BV-Moment}.
	
	(1)$\Leftrightarrow$(3) If $0\in \partial \|A(x)\|$ and $\lambda_{\max}(A(x))=-\lambda_{\min}(A(x))$, then using Proposition \ref{subdif norm and eig compact op},
	$$
	0\in {\rm co}\left(\partial \lambda_{\max}(A(x))\cup\partial \lambda_n(A(x)) \right)={\rm co}\left(m_{S_{\max}}\cup -m_{S_{\min}} \right)
	$$ 
	and with the same steps used on the proof of the (1)$\Leftrightarrow$(3) equivalence in Theorem \ref{teo equiv subdif, momentos, etc para matrices} follows that $m_{S_{\max}}\cap m_{S_{\min}}\neq \emptyset$. The converse can be proved similarly.
	
	To prove (2)$\Leftrightarrow$(3) we can use the formulas $\partial\lambda_{\max}(A(x))= m_{S_{\max}}$
	and $\partial \lambda_{\min}(A(x))=-m_{S_{\min}}$ proved in Theorem \ref{teo subdif y moment compact} and Proposition \ref{subdif norm and eig compact op}.
		Then it is trivial that   $m_{S_{\max}}\cap m_{S_{\min}}\neq \emptyset$ if and only if $0\in m_{S_{\max}}-m_{S_{\min}}=\partial \lambda_{\max}(A(x))+\partial \lambda_{\min}(A(x))$.
\end{proof}
\begin{remark}
	Observe that item (5) of Theorem \ref{equiv subdif, momentos, etc para compactos} allows the use of joint numerical ranges of finite self-adjoint matrices to decide the minimality of the compact operator.
\end{remark}
\section{Minimizing diagonals for rank one self-adjoint operators} \label{sec rank one}

Any rank-one self-adjoint (compact) operator $R\in B(H)^{sa}$ is a positive scalar multiple of an orthogonal projection $h h^{*}\in B(H)$ with  $h\in H$ and $\|h\|=1$. Then,  $D_0$ is a minimizing diagonal of $hh^{*}$ if and only if $r D_0$ is a minimizing diagonal of $r hh^{*}=R$. In this subsection we 
	will describe explicitly diagonals $D_0\in B(H)$ (in a fixed orthonormal basis $E=\{e_j\}_{j\in J}$ of $H$) such that $\|h h^{*}+D_0\|\leq \|h h^{*}+D\|$, for every diagonal (with respect to $E$) $D\in B(H)$. We will call them minimizing diagonals of $hh^*$ in the $E$ basis.
	We can suppose that $|h_j|> 0$, $\forall j$ and for numerable $j$, since otherwise we can work in a closed subspace of $H$. 
	In this context $(h_1,h_2,\dots,h_n,\dots)\in H$, with $h_j$ denotes the coordinates of $h$ in the fixed orthogonal basis $E$.
	 As mentioned, the results can be easily translated to general rank-one self-adjoint operators.

Let us consider $h\in H$ with $\|h\|=1$ and the rank one projection $h h^{*}\in B(H)$. 
We will explicitly describe diagonals $D_0\in B(H)$ (in a fixed orthonormal basis $E=\{e_j\}_{j\in J}$ of $H$) such that $\|h h^{*}+D_0\|\leq \|h h^{*}+D\|$, for every diagonal (with respect to $E$) $D\in B(H)$. We will call them minimizing diagonals of $hh^*$ in the $E$ basis.
We can suppose that $|h_j|> 0$, $\forall j$ and for numerable $j$, since otherwise we can work in a closed subspace of $H$. 
In this context $(h_1,h_2,\dots,h_n,\dots)\in H$, with $h_j$ denotes the coordinates of $h$ in the fixed orthogonal basis $E$.

%
%

The following is a slight generalization of the sufficient part of \cite[Theorem 2.2]{amlmrv} and its proof follows the same idea.
\begin{lemma}\label{lema Z minimal}
	Let $\mathcal{A}$ be a C$^*$-algebra, $\mathcal{B}\subset \mathcal{A}$ a C$^*$-subalgebra, $H$ a Hilbert space and $\rho:\mathcal{A}\to B(H)$ a representation of $\mathcal{A}$, and there exists $\xi\in H$, $\|\xi\|_H^2=1$, $Z\in \mathcal{A}$ such that $\langle \rho(Z)\xi,\rho(D)\xi\rangle=0$ $\forall D\in\mathcal{B}$, $\rho(Z^*Z)\xi=\|Z\|^2\xi$ then
	$$\|Z\|\leq \|Z+D\|\ , \forall D\in\mathcal{B}.
	$$
	That is, $Z$ is a minimal element with respect to $\mathcal{B}$.
\end{lemma}
\begin{proof}
	Observe that for every $D\in\mathcal{B}$
	\begin{equation*}
		\begin{split}
			\|Z+D\|^2&\geq \langle \rho(Z+D)\xi,\rho(Z+D)\xi\rangle\\
			&
			=\langle \rho(Z )\xi,\rho(Z )\xi\rangle+\langle \rho(Z )\xi,\rho( D)\xi\rangle+\langle \rho( D)\xi,\rho(Z )\xi\rangle+\langle \rho( D)\xi,\rho( D)\xi\rangle\\
			&\geq \langle \rho(Z )\xi,\rho(Z )\xi\rangle=\langle\rho(Z )^* \rho(Z )\xi,\xi\rangle=\langle\rho(Z^*   Z )\xi,\xi\rangle\\
			&=\|Z\|^2\langle\xi,\xi\rangle=\|Z\|^2
		\end{split}
	\end{equation*}
	and therefore $\|Z\|\leq \|Z+D\|\forall D\in\mathcal{B}$.
\end{proof}
We include here a result adapted to our needs.
\begin{lemma}\label{lema Z minimal con xi Hilbert-Schmidt}
	Let $Z$ be an operator of $B(H)$, and $\xi\in B_2(H)$ (a Hilbert-Schmidt operator) with $\tr(\xi^*\xi)=1$ such that $Z^*Z\xi=\|Z\|^2\xi$, $\tr(Z\xi(D\xi)^*)=\tr(Z\xi\xi^*D^*)=0$, $\forall D\in \text{Diag}(B(H))$ (the algebra of diagonal operators in a fixed basis), then
	$$\|Z\|\leq\|Z+D\|,\ \forall D\in \text{Diag}(B(H)).$$
\end{lemma}
\begin{proof}
	The proof is also motivated in the previous lemma.
	\begin{equation*}
		\begin{split}
			\|Z+D\|&\geq \tr\left((Z+D)\xi\left((Z+D)\xi\right)^*\right)=\tr(Z^*Z\xi\xi^*+ Z\xi\xi^* D^*+ D\xi\xi^* Z^*+D^*D\xi\xi^*)\\
			&=\tr(Z^*Z\xi\xi^*)+0+0+\tr(D^2\xi\xi^*)\geq \tr(Z^*Z\xi\xi^*)=\tr(\|Z\|^2\xi\xi^*)=\|Z\|^2\tr(\xi\xi^*)=\|Z\|^2
		\end{split}	
	\end{equation*}
	for all $ D\in \text{Diag}(B(H)).$
\end{proof}

The next result follows directly from \cite[Theorem 9]{KV 3x3} and \cite[Theorem 2]{bot var min length curves in unitary orb}. We state it here for the sake of clarity.
\begin{theorem}\label{teo generaliz de minimales de col gde}
	Let $T\in B(H)^{sa}$ described as an infinite matrix by $\left(T_{ij}\right)_{i,j\in \N}$ in a fixed basis.
	Suppose that $T$ satisfies that
	\begin{enumerate}
		\item[a) ] there exists $j_0\in \N$ satisfying $T_{j_0, j_0}=0$, with $T_{j_0, n}\neq 0$, for all $n\neq j_0$,
		\item[b) ] if $T^{(j_0)}$ is the operator $T$ with zero in its $j_0$th-column and $j_0$th-row then $$
		\left\|\col_{j_0}(T)\right\|\geq \left\|T^{(j_0)}\right\|
		$$
		(where $\left\|\col_{j_0}(T)\right\|$ denotes the Hilbert norm of the $j_0$th-column of $T$), and
		\item[c) ] $\left\langle \col_{j_0}(T),c_n(T)\right\rangle=0$ for each $n\in\N$, $n\neq j_0$.
	\end{enumerate}
	Then,
	\begin{enumerate}
		\item $\left\|T\right\|=\left\|\col_{j_0}(T)\right\|$.
		\item $T$ is minimal, that is
		$$
		\left\|T\right\|=\inf_{D\in \Diag(B(H)^{sa})}\left\|T+D\right\|=\inf_{D\in \Diag(K(H))}\left\|T+D\right\|,
		$$
		and $D=\Diag\big(\{T_{nn}\}_{n\in\N}\big)$ is the unique bounded minimal diagonal operator for $T$. 
	\end{enumerate}
\end{theorem}

Next, we introduce equivalent conditions for a rank one orthogonal projector in $B(H)$ to achieve minimality.

\begin{theorem} \label{teo minimal rango1}
	Let $h$ be an element of $H$ with $\|h\|_2=1$ and $h=(h_1,h_2,\dots,h_n,\dots)$ in a fixed basis $E$ of $H$. Then,
	\begin{enumerate}
		\item[(1)] if there exists $j_0$ such that $|h_{j_0}|^2>1/2$ then 
		$$hh^*- \Diag(1-|h_{j_0}|^2,\dots,1-|h_{j_0}|^2,|h_{j_0}|^2, 1-|h_{j_0}|^2,\dots)=$$ 
		$$=hh^*+(|h_{j_0}|^2-1)I+(1-2|h_{j_0}|^2) \ e_{j_0}e_{j_0}^*
		$$ 
		is a minimal matrix and is unique if $h_j\neq 0$ $\forall j$.
		\item[(2)] and if $|h_j|^2\leq 1/{2}$ for every $j$ then $D_0=-\frac12 I$ is a minimizing diagonal for $hh^*$. Moreover, if $h_j\neq 0$ $\forall j$, then this minimizing diagonal is unique (see also Corollary \ref{coro para hj0 nulo hay 2 diags minimales}).
	\end{enumerate}
\end{theorem}
\begin{proof}
	Recall that the diagonal of $hh^*\in H$ is $\Diag(hh^*)=\{|h_1|^2,|h_2|^2,\dots , |h_n|^2,\dots\}$ and $hh^*$ is a trace class positive operator (a projection or rank one) with $\tr(hh^*)=\sum_{j\in\N} |h_j|^2=1$ and hence a Hilbert-Schmidt operator with $\|hh^*\|_2=\tr(hh^*(hh^*)^*)=\tr(hh^*)=1$.
	We would also consider that the indexes $j$ belong to $\mathbb{N}$ although they could be finite in which case the proof is similar. We would also suppose that the coordinates $h_j\neq 0$ for all $j\in\N$ since otherwise the entire $j$-th row and column of $hh^*$ must be null and we can reorder the basis and take those $j$ away.
	\begin{enumerate}
		\item[(1)] 
		We will use Theorem \ref{teo generaliz de minimales de col gde} to prove that under the hypothesis $|h_{j_0}|^2>1/2$ the infinite matrix 
		$$
		m=hh^*+(|h_{j_0}|^2-1)I+(1-2|h_{j_0}|^2) \ e_{j_0}e_{j_0}^*
		$$
		is minimal. The diagonal of $m$ is 	 
		$$
		\Diag(m)=(|h_1|^2+|h_{j_0}|^2-1\ ,\ |h_2|^2+|h_{j_0}|^2-1, \dots ,\ \overbrace{0}^{j_0}\ ,|h_{j_0+1}|^2+|h_{j_0}|^2-1,\dots).
		$$
		Observe first that if $k\neq j_0$, since $\sum_{j} |h_j|^2=1$, then $m_{k,k}=|h_k|^2+|h_{j_0}|^2-1=-(1-|h_k|^2-|h_{j_0}|^2)=-\sum_{j\neq k, j_0} |h_j|^2$ and that  $m_{j_0,j_0}=|h_{j_0}|^2+(|h_{j_0}|^2-1)+(1-2|h_{j_0}|^2)=0$. With these elements in the diagonal a direct computation shows that the columns  $\col_k(m)$ and $\col_{j_0}(m)$ are orthogonal for $k\neq j_0$ (the elements of the diagonal were chosen for this purpose).
		
		Now consider the rank one operator $p_h^{(j_0)}=h^{(j_0)}\big(h^{(j_0)}\big)^*$, where $h^{(j_0)}$ equals $h$ except in the $j_0$ entry where there is a zero. Then its spectrum is $\sigma\left(p_h^{(j_0)}\right)=\{0,\|h^{(j_0)}\|^2\}=\{0,1-|h_{j_0}|^2\}$ and hence using functional calculus  $\sigma\left(p_h^{(j_0)}+(|h_{j_0}|^2-1)I^{(j_0)}\right)=\{|h_{j_0}|^2-1,0\}$, where $I^{(j_0)}$ is the identity matrix with a $0$ in the $j_0,j_0$ entry. Hence, denoting with $\col_{j_0}(m)$ the $j_0$-column of $m$, we have that
		\begin{equation}\label{eq norma mj0 menor que norma con cj0m}
			\begin{split}
				\|m^{(j_0)}\|&=\left\|p_h^{(j_0)}+(|h_{j_0}|^2-1)I^{(j_0)}\right\|=1-|h_{j_0}|^2=\sqrt{1-|h_{j_0}|^2}\sqrt{1-|h_{j_0}|^2}	\\
				&\leq \sqrt{1-|h_{j_0}|^2} |h_{j_0}|=\|\col_{j_0}(m)\|
			\end{split}
		\end{equation}
		where we used that $\sqrt{1-|h_{j_0}|^2}< |h_{j_0}|$ $\Leftrightarrow$ $1/2 < |h_{j_0}|^2$ and that $\col_{j_0}(m)=\overline{h_{j_0}}h^{(j_0)}$ then $\|\col_{j_0}(m)\|=|h_{j_0}| \|h^{(j_0)}\|=|h_{j_0}| \sqrt{\sum_{j\neq j_0} |h_j|^2}=\sqrt{1-|h_{j_0}|^2}$.
		
		Therefore, considering  that $\col_k(m)\perp\col_{j_0}(m)$ for $k\neq j_0$ and that $\|m^{(j_0)}\|\leq \|\col_{j_0}(m)\|$ (see \eqref{eq norma mj0 menor que norma con cj0m}) hold and Theorem \ref{teo generaliz de minimales de col gde} we can conclude that $m$ is a minimal matrix. Hence $\Diag(1-|h_{j_0}|^2,\dots,1-|h_{j_0}|^2,|h_{j_0}|^2, 1-|h_{j_0}|^2,\dots)$ is the closest diagonal to $hh^*$ if $h_j\neq 0$ for all $j$.

		\item[(2)] This item could be proved using Proposition \ref{prop K mas norma K por I mmal} but we include here a proof using other techniques regarding this special case. First we will show that there exists an element $k\in\left(\text{span}\{h\}\right)^\perp\subset H$ such that $|h_j|=|k_j|$ $\forall j\in \mathbb{N}$. This can be done considering an infinite polygon in the $\C$ plane with sides $|h_j|^2$ that starts and ends in $(0,0)$. This can be constructed if and only if $|h_j|^2\leq 1/{2}$, $\forall j\in\N$. Then define a collection of  angles $-\pi/2<\theta_j\leq\pi/2$, for $j\in \N$, of the corresponding sides of length $|h_j|^2$ with respect to the positive real axis required in order to obtain the mentioned closed polygon. 
		With these notations we obtain that 
		$$
		\sum_{j\in\N} e^{i\theta_j} |h_j|^2=0
		$$
		since the origin is where the polygon ends.
		Now, if $h_j=|h_j|e^{i\alpha_j}$, for $j\in \N$, we have that
		\begin{equation}
			\label{eq condicion que cumple hj poligono}		0=\sum_{j\in\N} e^{i\theta_j} |h_j|^2
			=\sum_{j\in\N} e^{i\alpha_j} |h_j| e^{i(\theta_j-\alpha_j)} |h_j|=\langle h,k\rangle
		\end{equation}
		
		for $k=\sum_{j\in\N}|h_j| e^{-i(\theta_j-\alpha_j)} e_j$. Hence $k \in \left(\text{span}\{h\}\right)^\perp$, satisfies that $|k_j|=|h_j|$ $\forall j\in \N$ and hence $\|k\|^2=\sum_{j\in\N}|k_j|^2=\sum_{j\in\N}|h_j|^2=1$.
		
		Then $hh^*$ is a rank one projector with eigenvectors $h$ and $k$ with corresponding eigenvalues $1$ and $0$  ($hh^*h=h$ and $hh^*k=h\langle k,h\rangle=0$).  Then the operator $Z=hh^*-\frac12 I$ has eigenvalues $\frac12$ and $\frac{-1}{2}$ with corresponding eigenvectors $h$ and $k$ (where $I\in B(H)$ denotes the identity operator). 
		
		Now consider the operator $Z=hh^*-\frac12 I$. Observe that the diagonal of $Z$ is $\big(|h_1|^2-1/2,|h_2|^2-1/2,\dots,|h_n|^2-1/2,\dots\big)$ in the fixed basis and that if we choose $\xi=\frac1{\sqrt{2}}(hh^*+kk^*)$ then $Z\xi=\frac1{2\sqrt{2}} (hh^*-kk^*)$ and $Z^*Z\xi=ZZ\xi=(\frac1{2})^2 \frac1{\sqrt2}(hh^*+kk^*)=\|Z\|^2\xi$. Moreover, using that $Z\xi^2=\frac14(hh^*-kk^*)$ and that the diagonal of $hh^*-kk^*$ is null, follows that  $\tr(Z\xi(D\xi)^*)=\tr(Z\xi^2D^*)=\frac1{4}\tr((hh^*-kk^*)D^*)=0$.
		Now we can apply Lemma \ref{lema Z minimal con xi Hilbert-Schmidt} with our defined $Z=hh^*-\frac12 I$ and $\xi=\frac1{\sqrt{2}}(hh^*+kk^*)$ to prove that $Z$ is a minimal operator with respect to $Diag(B(H))$.

	\end{enumerate}
\end{proof}
\begin{remark} Note that the minimizing diagonals for a rank one operator, as stated in Theorem \ref{teo minimal rango1}, are bounded but not compact.
\end{remark}

Next we show that the uniqueness of the minimizing diagonal fails if $h$ has any zero coordinate.

\begin{corollary}\label{coro para hj0 nulo hay 2 diags minimales}
	Let $h\in H$ such that $|h_j|\leq 1/2$ for all $j\in \N$, and suppose that  there exists $j_0\in \N$ such that $h_{j_0}=0$. Then, $hh^*\pm\frac 12e_{j_0}e_{j_0}^*$ are minimal operators. 
\end{corollary}
\begin{proof}
	By item 2 of Theorem \ref{teo minimal rango1}, $-\frac 12I$ is a minimizing diagonal for $hh^*$ and $\|hh^*-\frac 12I\|=\frac 12 $. Now consider $hh^*\pm\frac 12e_{j_0}e_{j_0}^*$, with $h$ as in the hypothesis. Then,
	$$\frac 12=\left\| \left( hh^*\pm\frac 12e_{j_0}e_{j_0}^*\right)e_{j_0} \right\| =\left\| C_{j_0}\left( hh^*\pm\frac 12e_{j_0}e_{j_0}^*\right) \right\| \leq \left\| hh^*\pm\frac 12e_{j_0}e_{j_0}^*\right\| $$
	and for each $j\neq j_0$,
	$$\left\| C_{j}\left( hh^*\pm\frac 12e_{j_0}e_{j_0}^*\right) \right\| =\sqrt{|h_j|^4+\sum_{k\neq j}|h_j|^2|h_k|^2}=\sqrt{|h_j|^4+|h_j|^2(1-|h_j|^2)}=|h_j|\leq \frac 12.$$
	Also, observe that 
	$$C_{j_0}\left( hh^*\pm\frac 12e_{j_0}e_{j_0}^*\right) \perp C_{j}\left( hh^*\pm\frac 12e_{j_0}e_{j_0}^*\right) $$
	for every $j\neq j_0$. Then, by Corollary 6.3 in \cite{BV Dga}, 
	$$\left\| C_{j_0}(hh^*\pm\frac 12e_{j_0}e_{j_0}^*)\right\| = \left\| hh^*\pm\frac 12e_{j_0}e_{j_0}^*\right\| =\frac 12=\left\| hh^*- \frac 12I\right\| ,$$
	therefore $hh^*\pm\frac 12e_{j_0}e_{j_0}^*$ are minimal operators. 
\end{proof}
The following is a related result with a different approach that provides conditions under which a diagonal matrix $D$ is minimal related to a rank-one operator.
\begin{lemma}
	Let $h\in H$ such that $\|h\|=1$ and $D\in D(B(H)^{sa})$. If there exists $j_0\in \N$ such that
	\begin{enumerate}
		\item $\|hh^*-D\|=\|(hh^*-D)e_{j_0}\|=\|C_{j_0}(hh^*-D)\|$.
		\item $|h_{j_0}|^2=D_{j_0,j_0}$.
	\end{enumerate}
	Then, $hh^*-D$ is minimal with $C_j(hh^*-D)\perp C_{j_0}(hh^*-D)$ for every $j\neq j_0$.
	
	Moreover, if $h_j\neq 0$ for all $j\in \N$, then $D$ is the unique minimizing diagonal and its entries are defined as
	
	\begin{equation}
		D_{jj}=|h_j|^2-\overline{h_{j}}h_{j_0}(1-|h_j|^2), \text{ for every } j\neq j_0.
	\end{equation}
\end{lemma}

\begin{proof}
	The minimality of $hh^*-D$ is a direct consequence of Lemma 6.1 in \cite{BV Dga}, since 
	$$(hh^*-D)_{j_0,j_0}=|h_{j_0}|^2-D_{j_0,j_0}=0.$$
	Moreover, if  
	$c_{j_0}(hh^*-D)_{j}=(hh^*-D)_{j,j_0}\neq0$ for 
	all $j\neq 
	j_0$, then $hh^*-D$
	has a unique 
	minimizing diagonal defined by
	\begin{equation*}
		(hh^*-D)_{j,j}=-\frac{\left\langle c_j(hh^*-D)_{\widecheck
				j},c_{j_0}(hh^*-D)_{\widecheck j}\right\rangle }{(hh^*-D)_{j,j_0} }, 
		\text{ for } j\neq j_0,
	\end{equation*}
	\\
	where $c_k(X)_{\widecheck l}\in 
	H\ominus 
	\text{gen}\{e_l\}$ is the  
	element obtained  after taking off 
	the $l^{\text{th}}$ entry of 
	$c_k(X)\in H$. Then, for $j\neq j_0$
	\begin{eqnarray*}
		(hh^*-D)_{j,j}&=&-\frac{\left\langle c_j(hh^*-D)_{\widecheck
				j},c_{j_0}(hh^*-D)_{\widecheck j}\right\rangle }{(hh^*-D)_{j,j_0} }\\
		|h_j|^2-D_{jj}&=&\sum_{i\neq j}|h_i|^2\overline{h_{j}}h_{j_0}+\overline{h_{j}}h_{j_0}\left( |h_{i_0}|^2-D_{i_0,i_0}\right) \\
		|h_j|^2-D_{jj}&=&\sum_{i\neq j}|h_i|^2\overline{h_{j}}h_{j_0}\\
		D_{jj}&=&|h_j|^2-\sum_{i\neq j}|h_i|^2\overline{h_{j}}h_{j_0}\\
		D_{jj}&=&|h_j|^2-\overline{h_{j}}h_{j_0}(1-|h_j|^2)\\
	\end{eqnarray*}
\end{proof}

\end{document}